\documentclass{article}

\usepackage{amssymb}
\usepackage{graphicx}
\usepackage{setspace}
\usepackage{amsmath}
\usepackage{float}
\usepackage{color}
\usepackage{graphicx}
\usepackage{float}
\usepackage{algorithm} 
\usepackage{algorithmic}  
\usepackage[algo2e]{algorithm2e} 
\providecommand{\keywords}[1]{\textbf{\textit{Keywords:}} #1}
\usepackage[left=2.5cm, right=2.5cm, top=2cm]{geometry}
\begin{document}



\title{An approximate It\^o-SDE based simulated annealing algorithm for multivariate design optimization problems}

\author{A. Batou$^{a \ast}$ \\  $^{a}$Universit\'e Paris-Est, Laboratoire Mod\'elisation et Simulation Multi Echelle,\\ MSME UMR 8208 CNRS, France}

\maketitle

\begin{abstract}
This research concerns design optimization problems involving numerous design parameters and large computational models. These problems generally consist in non-convex constrained optimization problems in large and sometimes complex search spaces. The classical simulated annealing algorithm rapidly loses its efficiency in high search space dimension. In this paper a variant of the classical simulated annealing algorithm is constructed by incorporating (1) an It\^o stochastic differential equation generator (ISDE) for the transition probability and (2) a polyharmonic splines interpolation of the cost function. The control points are selected iteratively during the research of the optimum. The proposed algorithm explores efficiently the design search space to find the global optimum of the cost function as the best control point. The algorithm is illustrated on two applications. The first application consists in a simple function in relatively high dimension. The second is related to a Finite Element model.\\

\keywords{Engineering Design; Simulated annealing; It\^o stochastic differential equation; polyharmonic spline interpolation)}

\end{abstract}

\section{Introduction}

Recent advances in computational mechanics offer engineers the possibility of constructing very predictive computational models of complex mechanical structures. Such computational models can then be used to assess the performance of the structure over specific design scenarios. This increase of predictability is inevitably accompanied by an increase of the size of the computational models. In most industrial sectors, it is common for engineers to deal with computational models having tens of millions of degrees of freedom. The forward problem consisting in calculating the response of an already designed system can be performed using usual computational resources. Nevertheless the design of a system represented by a large computational model is a hard task since it generally consists in a large non-convex constrained optimization problem in a large and sometimes complex search space and requires numerous runs of the cost function.

There exists several methods in the literature in order to find or approximate the global maximum of a non-convex constrained optimization problem \cite{Pardalos2002}. Among these methods, simulated annealing methods \cite{Kirkpatrick1983,Cerny1985} are attractive metaheuristic random search methods which consist in assigning the design parameters a probability distribution whose logarithm is proportional to the cost function to be minimized. Then a Markov Chain is generated and by varying the coefficient of proportionality (inverse of the temperature) in an appropriate way, the Markov Chain converges to the global optimum. This methodology, for which the basic version \cite{Kirkpatrick1983,Cerny1985} is easy to implement, has been improved by many researchers (see for instance \cite{Pardalos2002,Ingber1993,Suman2006}). In \cite{Zolfaghari1999} the authors have combined the simulated annealing with a tabu search approach in order to avoid solution re-visits. In \cite{Mingjun2004} the authors introduce a chaotic initialization and a chaotic sequence generation to improve the convergence rate of the simulated annealing method. Many researches have been devoted to the extension of simulated annealing method to multiobjective optimization problems \cite{Suman2006}. The improvement of the annealing schedule has also interested many researchers (see for instance \cite{Kouvelis1992} for the choice of the initial temperature  and \cite{Ingber1989,Noutani1998,Azizi2004,Triki2004}  for the  cooling schedule). Convergence results related to simulated annealing algorithm have been studied in \cite{Aarts1985,Lundy1986,Faigle1988,Faigle1991,Granville1994}  for the case where the temperature is held constant until a stationary distribution is reached  and in  \cite{Mitra1986,Hajek1988,Connors1989,Borkar1992,Anily1987,Yao1991} for the case where the temperature always changes (and then no stationary distribution can be reached). For cost functions involving large computational models, the simulated annealing algorithm can become prohibitively time consuming. To circumvent this difficulty, an adaptative surface response can be introduced for the approximation of the cost function \cite{Wang2001}.

The basic version of the simulated annealing \cite{Kirkpatrick1983,Cerny1985} is based on a Metropolis-Hastings algorithm \cite{Hastings1970} for which a proposal distribution has to be provided. If the dimension of the search spaces increases, it is not easy to defined this proposal distribution and the Markov chain can become very slow due to numerous rejections. In this case, for large computational models, the use of the classical simulated annealing algorithm can become prohibitive due to the large number of runs of the cost function that are needed.

The objective of this paper is to improve the simulated annealing algorithm in order to design, in high-dimensional search space, complex mechanical systems represented by large computational models. The methodology proposed in the present paper is based on two ingredients. The first one consists in replacing the Metropolis-Hastings algorithm by an It\^o stochastic differential equation (ISDE) generator \cite{Soize1994,Soize2008} which belongs to class of Monte Carlo Markov Chain Methods (MCMC). This generator is very well adapted in high dimension and does not require the specification of a proposal distribution. Nevertheless it requires the gradient of the cost function to be evaluated thousands times making it unusable. To circumvent this difficulty, the second ingredient of the methodology proposed here consists, as in \cite{Wang2001}, in introducing response surface function in order to approximate the cost function. Since we are interested in multivariate optimization problems, a radial interpolation surface (polyharmonic) is preferred to a polynomial approximation (as used in\cite{Wang2001}). This interpolation is well adapted in relatively high input dimension and allows the gradient of the cost function to be calculated explicitly at each step of the Markov Chain. An adaptative strategy is proposed in order to enrich sequentially, during the numerical integration of the ISDE, the set of the control points that are used to construct the interpolation surface. At the end of the simulation the best control point is selected as the global optimum. 

In Section \ref{sec2}, the design optimization problem and the classical simulated annealing algorithm are presented. Then Section \ref{sec3} is devoted to the presentation the ISDE-based simulated annealing algorithm. The constraints on the design parameters are taken into account by introducing a regularisation of the indicator function. In Section \ref{sec4}, the polyharmonic interpolation of the cost function is introduced and the complete approximate ISDE-based simulated annealing algorithmdesign optimization algorithm is presented. The methodology is illustrated through two applications in Section \ref{sec5}. The first application consists in the classical Ackley's test function in high dimension. The second one consists in a Finite Element structure submitted to an acceleration of the soil and for which the maximal displacement has to be minimized by designing the stiffnesses of supporting springs.

\section{Design optimization problem and classical simulated annealing algorithm}
\label{sec2}
In this section the design optimization problem is formulated and the simulated annealing probability distribution is constructed classically using the Maximum Entropy (MaxEnt) principle.

Let $\boldsymbol{a} $ be the vector of the design parameters with values in $\mathbb{R}^{N}$.
Let $\mathcal{C}_a \subseteq \mathbb{R}^{N}$ be the admissible set for vector $\boldsymbol{a}$. 
It is assumed that the design optimization problem to be solved consists in searching the global minimum of a non-convex positive cost function $\mathcal{D}(\boldsymbol{a})$. Then the optimal value $\boldsymbol{a}^{\rm opt}$ of vector $\boldsymbol{a}$ is calculated such that
\begin{eqnarray}
\label{eq1}
\boldsymbol{a} ^{\textrm{opt}} = \arg  \min_{ \boldsymbol{a}\in \mathcal{C}_a} \mathcal{D}(\boldsymbol{a}) \, .
\end{eqnarray}
For instance,  in case of a least square problem, the function  $\mathcal{D}(\boldsymbol{a})$ can be written as
\begin{eqnarray}
\label{eq2}
\mathcal{D}(\boldsymbol{a}) = \Vert \textbf{G}(\boldsymbol{a}) - \textbf{g}^* \Vert^2   \, ,
\end{eqnarray}
in which $\textbf{G}(\boldsymbol{a})$ is a performance vector-valued function and $\textbf{g}^*$ is a target vector. 
In random search methods, a random vector $\boldsymbol{A}$ with values in $\mathcal{C}_a$ is associated with the vector $\boldsymbol{a}$. In the simulated annealing the probability density function (pdf) $p_{\boldsymbol{A}}(\boldsymbol{a})$ of random vector $\boldsymbol{A}$ is constructed by using the MaxEnt principle \cite{Shannon1948,Jaynes1957,Kapur1992} under the constraints defined by the available information related to random vector $\boldsymbol{A}$.
This information is written as
\begin{eqnarray}
\label{eq3}
\begin{split}
\boldsymbol{a} \in \mathcal{C}_a \, , \quad \quad (a)\\
\int_{\mathcal{C}_a} {\mbox{1\hspace{-.25em}l}_{ \mathcal{C}_a } (\boldsymbol{a}) \, \exp({-\lambda\, \mathcal{D}(\boldsymbol{a})})}\,d\boldsymbol{a}) = 1 \, ,  \quad \quad (b)\\
E \{ \mathcal{D}(\boldsymbol{a}) \} = c < +\infty \,, \quad \quad (c)
\end{split}
\end{eqnarray}
in which $E \{ . \}$ is the mathematical expectation. Equation(\ref{eq3}-b) is related to the normalization of the pdf and Eq.~(\ref{eq3}-c) is related to the finiteness of the mean value of the cost function which is interpreted as an energy in thermodynamics. The entropy related to pdf $p_{\boldsymbol{A}}$ is defined by
\begin{eqnarray}
\label{eq4}
S(p_{\boldsymbol{A}}) = -\int_{\mathbb{R}^{N}}{ p_{\boldsymbol{A}}(\boldsymbol{a}) \, \textrm{log} ( p_{\boldsymbol{A}}(\boldsymbol{a})) \, \textrm{d} \boldsymbol{a}} \, ,
\end{eqnarray}
where $\log$ is the natural logarithm. This functional measures the relative uncertainty of $p_{\boldsymbol{A}}$. Let $\mathcal{C}$ be the set of all the possible pdfs of random vector $\boldsymbol{A}$, satisfying the constraints defined by Eq.~(\ref{eq3}). Then the MaxEnt principle consists in constructing the pdf $p_{\boldsymbol{A}}$ as the unique pdf in $\mathcal{C}$ which maximizes entropy $S(p_{\boldsymbol{A}})$. Then by introducing a positive Lagrange multiplier $\lambda_0$ associated with Eq.~(\ref{eq3}-b) and a Lagrange multiplier ${\lambda}$  associated with Eq.~(\ref{eq3}-c), it can be shown (see \cite{Jaynes1957,Kapur1992}) that the MaxEnt solution, is defined by
\begin{eqnarray}
\label{eq5}
p_{\boldsymbol{A}}(\boldsymbol{a})  = \mbox{1\hspace{-.25em}l}_{ \mathcal{C}_a } (\boldsymbol{a}) \, \lambda_0 \, \exp({-\lambda\, \mathcal{D}(\boldsymbol{a})})
\end{eqnarray}
in which $\lambda_0 = (\int_{\mathcal{C}_a} {\mbox{1\hspace{-.25em}l}_{ \mathcal{C}_a } (\boldsymbol{a}) \, \exp({-\lambda\, \mathcal{D}(\boldsymbol{a})})}\,d\boldsymbol{a})^{-1}$ is a normalization constant. In the context of the simulated annealing algorithm, the Lagrange multiplier is written as $\lambda = \frac{1}{T}$ (or $\lambda = \frac{1}{\beta T}$ with $\beta>0$) in which $T$ is the temperature which controls the sharpness of pdf $p_{\boldsymbol{A}}$. Then random variable $\boldsymbol{A}$ is rewritten $\boldsymbol{A}_T$ and its pdf is rewritten as 
\begin{eqnarray}
\label{eq6}
p_{\boldsymbol{A}_T}(\boldsymbol{a})  = \mbox{1\hspace{-.25em}l}_{ \mathcal{C}_a } (\boldsymbol{a}) \, \lambda_0 \, \exp({-\frac{ \mathcal{D}(\boldsymbol{a})}{T}})
\end{eqnarray}
A large value of the temperature makes $p_{\boldsymbol{A}_T}$ tend to the uniform distribution and a small value of the temperature makes $p_{\boldsymbol{A}_T}$ tend to a Dirac distribution centred at the optimal value. Then the classical simulated annealing algorithm consists in constructing a sequence $(\boldsymbol{A}^{1}, \boldsymbol{A}^{2}, \ldots, \boldsymbol{A}^{{n_t}})$ using the Metropolis-Hastings (MH) algorithms \cite{Hastings1970} with the pdf $p_{\boldsymbol{A}_T}$ while decreasing slowly the temperature (see for instance \cite{Noutani1998} for the cooling strategies). The algorithm is summarized in Algorithm \ref{algo1}. 
\begin{algorithm}[h!]
\textbf{INITIALIZATION:}\\
Generate initial value $\boldsymbol{A}_{0}$\;
$\boldsymbol{A}^{1} = \boldsymbol{A}_{0}$\;
\textbf{LOOP:}\\
	\For{$k=1,\ldots,(n_t-1)$}{
		Set temperature $T_k$\,\;\\
		Set the size $M_k$ of the MH sequence \,\;\\
		MH algorithm: \\
	\For{$\ell=1,\ldots,(M_k-1)$}{
		Generate a neighbour $ \boldsymbol{A}^{k,new}$ of $\boldsymbol{A}^{k}$\,\;\\
		\If{$\mathcal{D}( \boldsymbol{A}^{k,new}) < \mathcal{D}( \boldsymbol{A}^{k})$}{
		$\boldsymbol{A}^{k} \leftarrow \boldsymbol{A}^{k,new}$\,\;
		}
		\Else{
		Generate $u$ uniformly in $[0,1]$ \,\;\\
		\If{$u<\exp({-{ (\mathcal{D}( \boldsymbol{A}^{k,new}) - \mathcal{D}(\boldsymbol{A}^{k}))}/{T_k}})$}{
		$\boldsymbol{A}^{k} \leftarrow  \boldsymbol{A}^{k,new}$ \,\;
		}
		}
		}
		$\boldsymbol{A}^{k+1} \leftarrow \boldsymbol{A}^{k}$\,\;

}
\caption{Classical simulated annealing algorithm }
\label{algo1}
\end{algorithm}
In this algorithm \ref{algo1}, it is often set $M_k = 1$ for $k=1\ldots,n_t$, skipping the burning phase. In this case, there is only one loop and the temperature is changed at each iteration.
Under some conditions related to the rate of decreasing of the temperature, the algorithm converges to the global optimum.
Two difficulties are usually referred for the classical simulated annealing algorithm:(1) it requires numerous calculation of the cost function $\mathcal{D}(\boldsymbol{a})$ making it unusable for large computational model and (2) if the dimension of random design vector $\boldsymbol{A}_T$ is large, the Metropolis-Hastings algorithms may have difficulties to explore efficiently the design search space. In the next sections two modifications are introduced to circumvent these two difficulties. The first one consists in replacing the Metropolis-Hastings generator by an ISDE generator in order to explore efficiently the design search space in high dimension. The second one consists in introducing a polyharmonic spline interpolation of the function $\boldsymbol{a} \mapsto \mathcal{D}(\boldsymbol{a})$ for which the list of control points is enriched during the algorithm.
\section{Construction of a new ISDE-based simulated annealing algorithm}
\label{sec3}
The ISDE generator consists in constructing the pdf of random vector $\boldsymbol{A}_T$ as the density of the invariant measure $p_{\boldsymbol{A}_T}(\boldsymbol{a})\textrm{d} \boldsymbol{a}$ associated with the stationary solution of a second-order nonlinear ISDE. This methodology, which has similarities with the Hamiltonian (or Hybrid) Monte Carlo method \cite{Duane1987,Salazar1997}, has recently been revisited \cite{Soize1994,Soize2008} in the context of the generation of random vectors in high dimension for which the pdf is constructed using the MaxEnt principle with complex available information. The advantages of this generator compared to the other MCMC generators such as the Metropolis-Hastings algorithm (see \cite{Hastings1970}) are: (1) there is no need to introduce a proposal distribution (uniform in case of random walk) for the exploration of the admissible space, (2) a damping can be introduced in order to rapidly reach the invariant measure and (3) mathematical results concerning ISDEs can be used for convergence analyses.
Below, we first present the ISDE generator followed by its application for generating independent realizations of $\boldsymbol{A}_T$. Then a new simulated algorithm based on the use of this ISDE generator is presented
\subsection{ISDE generator}
\label{sec31}
%


Let $\Psi(\boldsymbol{u})$ be a positive function. 
Let $\{(\boldsymbol{U}(r), \boldsymbol{V}(r)), r \geq 0 \}$ be a stochastic process  with values in $\mathbb{R}^N \times \mathbb{R}^N$ satisfying, for all $r \geq 0$, the following ISDE  \cite{Ito1951,Soize1994}
\begin{eqnarray}
\label{eq8}
\renewcommand{\arraystretch}{2.8}
\begin{array}{ccl}
\textrm{d} \boldsymbol{U}(r) & = & \boldsymbol{V}(r)\, \textrm{d}r \\[-.4cm]
\textrm{d} \boldsymbol{V}(r) & = & -\boldsymbol{\nabla}_{\boldsymbol{u}} \Psi(\boldsymbol{U})\, \textrm{d}r -\frac{1}{2} [D]\boldsymbol{V}(r)\, \textrm{d}r +[S]\, \textrm{d}\boldsymbol{W}(r) \, ,
\end{array}
\end{eqnarray}
with the initial conditions
\begin{eqnarray}
\label{eq9}
\boldsymbol{U}(0) = \boldsymbol{U}_{0} \quad , \quad \boldsymbol{V}(0) = \boldsymbol{V}_{0} \quad a.s \, .
\end{eqnarray}
In Eq.~(\ref{eq8}), $\boldsymbol{\nabla}_{\boldsymbol{u}} \Psi(\boldsymbol{u})$ is the gradient of function $\Psi(\boldsymbol{u})$ with respect to $\boldsymbol{u}$. The matrix $[D]$ is a symmetric positive-definite damping matrix, the lower triangular matrix $[S]$ is  such that $[D] = [S][S]^T$ and $\boldsymbol{W} = (W_1,\ldots,W_N)$ is the normalized Wiener stochastic process indexed by $\mathbb{R}^+$. The random initial condition $(\boldsymbol{U}_{0},\boldsymbol{V}_{0})$ is a second-order random variable independent of the Wiener stochastic process $\{\boldsymbol{W}(r) , r \geq 0\}$.
Then it can be proven (see \cite{Soize1994,Soize2008}) that, if $\boldsymbol{u}\mapsto \Psi(\boldsymbol{u})$ is continuous on $\mathbb{R}^N$, if
$\boldsymbol{u}\mapsto \Vert \boldsymbol{\nabla}_{\boldsymbol{u}}\Psi(\boldsymbol{u}) \Vert$ is locally bounded on $\mathbb{R}^N$
(i.e. is bounded on all compact set in $\mathbb{R}^N$), and if
\begin{eqnarray}
\label{eq9a}
\inf_{\Vert \boldsymbol{u} \Vert > R} \Psi(\boldsymbol{u}) \rightarrow +\infty \quad \textrm{if} \, R \rightarrow +\infty \, ,
\end{eqnarray}
\begin{eqnarray}
\label{eq9b}
\inf_{ \boldsymbol{u} \in \mathbb{R}^N} \Psi(\boldsymbol{u}) = \Phi_{\textrm{min}}  \quad \textrm{with} \, \Phi_{\textrm{min}} \in \mathbb{R} \, ,
\end{eqnarray}
\begin{eqnarray}
\label{eq9c}
\int_{\mathbb{R}^{N}} {\Vert \boldsymbol{\nabla}_{\boldsymbol{u}}\Psi(\boldsymbol{u}) \Vert \,  \exp(-\Psi(\boldsymbol{u}))\, \textrm{d} \boldsymbol{u} }\, < \,+\infty \, ,
\end{eqnarray}
then the ISDE defined by Eqs.~(\ref{eq8}) and (\ref{eq9}) admits an invariant measure defined by the pdf
 $\rho(\boldsymbol{u},\boldsymbol{v})$ which is written as
\begin{eqnarray}
\label{eq10}
\rho(\boldsymbol{u},\boldsymbol{v}) = \lambda_0\,\textrm{exp} (-{\boldsymbol{\nabla}_{\boldsymbol{u}} \Psi(\boldsymbol{u})} \times (2\pi)^{-N/2}\, \textrm{exp} (-\frac{1}{2} \Vert \boldsymbol{v} \Vert ^2) \, .
\end{eqnarray}
It can then be deduced that,  for $r \rightarrow +\infty$, the stochastic process $\{\boldsymbol{U}(r) , r \geq 0\}$ tends to a stationary stochastic process in probability distribution, for which the one-order marginal probability distribution is $p_{\boldsymbol{U}^{\rm st}}(\boldsymbol{u})$ such that
\begin{eqnarray}
\label{eq11}
p_{\boldsymbol{U}^{\rm st}}(\boldsymbol{u})  =  \lambda_0 \, \exp(-\Psi(\boldsymbol{u})) \,.
\end{eqnarray}
Therefore, using an independent realization of the Wiener stochastic process $\boldsymbol{W}$ and an independent realization of the initial condition $(\boldsymbol{U}_{0},\boldsymbol{V}_{0})$, an independent realization random variable $\boldsymbol{U}^{\rm st}$ related to the probability distribution (\ref{eq11}) can be constructed as the solution of the ISDE defined by Eqs.~(\ref{eq8}) and (\ref{eq9}), for $r$ sufficiently large. In other words, independent realizations of any random variable having a probability distribution of the form (\ref{eq11}) (verifying the required continuity and regularity properties for function $\Psi$) can be obtained by solving the ISDE defined by Eqs.~(\ref{eq8}) and (\ref{eq9}).
The value $r_0$ of $r$ for which the invariant measure is approximatively reached depends on the choice of the damping matrix $[D]$ and on the probability distribution of the random initial conditions. The damping induced by the matrix $[D]$ has to be sufficiently large in order to rapidly kill the transient response but a too large damping introduces increasing errors in the numerical integration of the ISDE. 

The solution of the ISDE can be constructed numerically using a for instance a modified Euler integration schemes which has a better accuracy than a classical Euler integration schemes (see \cite{Soize2008,Burrage2007}).
Let $\Delta r$ be the integration step size and let  $\{ r_\ell, \ell=1,\ldots , M\}$ be the sampling points, M being a positive integer. Let  $\{ {U}^{\ell}, \ell=1,\ldots , M\}$ and $\{ {V}^{\ell}, \ell=1,\ldots , M\}$ be the time series related to the discretization of random processes  $\boldsymbol{U}(r)$ and  $\boldsymbol{V}(r)$.
Then, after having generated the random initial values ${U}^{1} = \boldsymbol{U}_{0}$ and ${V}^{1} = \boldsymbol{V}_{0}$, the following updating rules applies
\begin{eqnarray}
\label{eq11b}
\renewcommand{\arraystretch}{2.8}
\begin{array}{ccl}
		\boldsymbol{V}^{\ell+1} &  = &  ([I]-\frac{\Delta r}{2} [D]) \boldsymbol{V}^{\ell}- \Delta r \boldsymbol{\nabla}_{\boldsymbol{u}} \Psi (\boldsymbol{U}^{\ell})  +[S] \Delta \boldsymbol{W}^{\ell+1}\,, \\[-.4cm]
		\boldsymbol{U}^{\ell+1} &  = &  \boldsymbol{U}^{\ell} + \Delta r\, \boldsymbol{V}^{\ell+1}\,\;
\end{array}
\end{eqnarray}
The vector $\Delta \boldsymbol{W}^{\ell+1}$ is a second-order Gaussian centered random vector with covariance matrix equal to $\Delta r \, [I_N]$ and the random vectors $\Delta \boldsymbol{W}^{1}, \ldots,\Delta \boldsymbol{W}^{M}$ are mutually independent.
\subsection{ISDE generator for random variable $\boldsymbol{A}_T$}
\label{sec32}
Equation (\ref{eq6}) can be rewritten as
\begin{eqnarray}
\label{eq12}
p_{\boldsymbol{A}_T}(\boldsymbol{a})  =  \, \lambda_0 \, \exp( \log( \mbox{1\hspace{-.25em}l}_{ \mathcal{C}_a } (\boldsymbol{a}))  {-\frac{ \mathcal{D}(\boldsymbol{a})}{T}})
\end{eqnarray}
Comparing equations (\ref{eq11}) and (\ref{eq12}), we naturally introduce the function
\begin{eqnarray}
\label{eq13}
\Psi_T(\boldsymbol{a}) =  {\frac{ \mathcal{D}(\boldsymbol{a})}{T}} - \log( \mbox{1\hspace{-.25em}l}_{ \mathcal{C}_a } (\boldsymbol{a})) \, ,
\end{eqnarray}
Unfortunately, due to the presence of the indicator function $\boldsymbol{a} \mapsto \mbox{1\hspace{-.25em}l}_{ \mathcal{C}_a } (\boldsymbol{a})$, the function $\Psi_T(\boldsymbol{u})$ is not differentiable and then cannot be used within an ISDE generator. To solve this issue, the indicator function $\boldsymbol{a} \mapsto \mbox{1\hspace{-.25em}l}_{ \mathcal{C}_a } (\boldsymbol{a})$  is replaced by the regularized indicator function  $\boldsymbol{a} \mapsto \mbox{1\hspace{-.25em}l}_{ \mathcal{C}_a } ^{\rm{reg}}(\boldsymbol{a})$. For simple manifolds, explicit regularized indicator function can be constructed. Two simple examples are given below (see \cite{Batou2014}):

{\parindent 0pt
\textbf{Example 1}: Positive design variables

Assume that the constraints on $\boldsymbol{a}$ are
\begin{eqnarray}
\label{eq14}
a_i > 0 , \quad \quad i=1\ldots,N   \,.
\end{eqnarray}
In this case, the regularized indicator function can be constructed as follow
\begin{eqnarray}
\label{eq15}
\mbox{1\hspace{-.25em}l}_{ \mathcal{C}_a } ^{\rm{reg}}(\boldsymbol{a}) = \prod_{i=1}^N \frac{1}{2}\left(1+\tanh\left(\frac{a_i}{\alpha_i}\right)\right)  \,,
\end{eqnarray}
where $\alpha_i >0, i=1\ldots,N$ are shape parameters 

\textbf{Example 2}: Design variables belonging to independent intervals

Assume that the constraints on $\boldsymbol{a}$ are
\begin{eqnarray}
\label{eq16}
l_i < a_i < u_i , \quad \quad i=1\ldots, N \,,
\end{eqnarray}
where $l_i$ and $u_i$, $i=1\ldots, N$ are lower and upper bounds. In this case, the regularized indicator function can be constructed as follow
\begin{eqnarray}
\label{eq17}
\mbox{1\hspace{-.25em}l}_{ \mathcal{C}_a } ^{\rm{reg}}(\boldsymbol{a}) = \prod_{i=1}^N \frac{1}{4}\left(1+\tanh\left(\frac{a_i - l_i}{\alpha_i}\right)\right)\left(1+\tanh\left(\frac{u_i - a_i}{\alpha_i}\right)\right)  \,.
\end{eqnarray}
For more complex manifolds, a general kernel-smoothing regularization approach has been proposed in \cite{Guilleminot2014}.

Then a regularized ISDE can be introduced by replacing, in Eq.(\ref{eq8}), $\Psi$ by $\Psi_T ^{\rm{reg}}$ such that
\begin{eqnarray}
\label{eq18}
\Psi_T ^{\rm{reg}}(\boldsymbol{a}) =  {\frac{ \mathcal{D}(\boldsymbol{a})}{T}} - \log( \mbox{1\hspace{-.25em}l}_{ \mathcal{C}_a } ^{\rm{reg}} (\boldsymbol{a})) \, ,
\end{eqnarray}
and then independent realizations of random variable $\boldsymbol{A}_T$ can be constructed.
\subsection{ISDE based simulated annealing algorithm}
\label{sec33}
Now that an ISDE generator for random variable $\boldsymbol{A}_T$ is constructed, we can use this generator for replacing the MH algorithm in the classical simulated annealing algorithm \ref{algo1}. This ISDE based simulated annealing algorithm is summarized in Algorithm \ref{algo2}. 

\begin{algorithm}[h!]
\textbf{INITIALIZATION:}\\
Generate $\boldsymbol{U}_{0}$ and $\boldsymbol{V}_{0}$\,\;\\
$\boldsymbol{A}^{1} = \boldsymbol{U}_{0}$\,\;\\
$\boldsymbol{U}^{1} = \boldsymbol{U}_{0}$\,\;\\
$\boldsymbol{V}^{1} = \boldsymbol{V}_{0}$\,\;\\
\textbf{LOOP:}\\
	\For{$k=1,\ldots,(n_t-1)$}{
	Set temperature $T_k$\,\;\\
	Set damping matrix $[D_k] = [S_k][S_k]^T$\,\;\\
	Set increment $\Delta r_k$\,\;\\
	Set number of steps $M_k$\,\;\\
	Solve ISDE : \\
	\For{$\ell=1,\ldots,(M_k-1)$}{
		$\boldsymbol{V}^{\ell+1}   =   ([I]-\frac{\Delta r_k}{2} [D_k]) \boldsymbol{V}^{\ell}- \Delta r_k \boldsymbol{\nabla}_{\boldsymbol{u}} \Psi_{T_k}^{\rm{reg}} (\boldsymbol{U}^{\ell})  +[S_k] \Delta \boldsymbol{W}^{\ell+1}$\,\;\\
		$\boldsymbol{U}^{\ell+1}   =   \boldsymbol{U}^{\ell} + \Delta r_k\, \boldsymbol{V}^{\ell+1}$\,\;\\
		}
		$\boldsymbol{A}^{k+1} = \boldsymbol{U}_{M_k}$\,\;\\
		$\boldsymbol{U}^{1} \leftarrow \boldsymbol{U}_{M_k}$\,\;\\
		$\boldsymbol{V}^{1} \leftarrow \boldsymbol{V}_{M_k}$\,\;\\
}
\caption{ISDE-based simulated annealing algorithm}
\label{algo2}
\end{algorithm}
At each changing of the temperature, three parameters have to be set: the steps in the ISDE $\Delta r_k$, the damping matrix $[D_k] = [S_k][S_k]^T$ and the number of steps $M_k$ in the ISDE. By analogy to second-order dynamical systems, the time step that guaranties the stability of the modified Euler scheme is calculated using the eigenvalues related the conservative part of the dynamical equation. In the ISDE (\ref{eq8}), the corresponding "mass" matrix is identity and the "stiffness" matrix for given "displacement" $\boldsymbol{U}$ is the hessian matrix $[H_{\Psi_T}(\boldsymbol{U})]$ of function $\Psi_T(\boldsymbol{U})$. Then at each iteration $k$, the step is calculated as $\Delta r_k = 2\,\pi/(m\sqrt{\lambda_{max}})$ where $m$ is an integer such that $m>10$ and $\lambda_{max}$ is the largest eigenvalue of the hessian matrix $[H_{\Psi_T}(\boldsymbol{U}^1)]$. When approaching the border of the admissible domain $\mathcal{C}_a$,  the ISDE can become very "stiff", which may yield stability issues. To avoid this, the integer $m$ has to be chosen large enough but not too large in order to avoid a too slow integration of the ISDE. Alternatively, an adaptative step-size algorithm can be used to automatically choose the best step-size at each iteration of the Euler Scheme \cite{Guilleminot2014}.
The damping matrix has to be large enough in order to kill rapidly the transient response. But a too large damping will introduce errors during the integration. A good compromise consists in imposing the largest damping rate equal to $\xi = 0.7$. This can achieved approximatively by choosing $[D_k] = 2\,\xi\sqrt{\lambda_{max}}[I_N]$. The number of steps of $M_k$ can be chosen as invariant and then determined through an off-line convergence analysis. The convergence can also be controlled during the integration of each ISDE yielding different values for $M_1, \ldots, M_k$.
\section{Approximation of the cost function}
\label{sec4}

In this section, the function $\boldsymbol{a} \mapsto \mathcal{D}(\boldsymbol{a})$ in the ISDE is approximated using multivariate polyharmonic splines. This approximation allows to reduce the computational cost for large computational models. Furthermore it allows to derive explicitly the gradient $\boldsymbol{a} \mapsto \boldsymbol{\nabla}_{\boldsymbol{a}} \mathcal{D}(\boldsymbol{a})$, avoiding the calculation of the gradient numerically. The multivariate polyharmonic splines approximation $\boldsymbol{a} \mapsto \mathcal{D}^{\rm mps}(\boldsymbol{a})$ of function $\boldsymbol{a} \mapsto \mathcal{D}(\boldsymbol{a})$ is written as
\begin{eqnarray}
\label{eq29}
\mathcal{D}^{\rm mps}(\boldsymbol{a}) = \boldsymbol{w}^T \boldsymbol{b}(\boldsymbol a) ,
\end{eqnarray}
in which $\boldsymbol{b}(\boldsymbol a) = (b_1(\boldsymbol a), \ldots b_{n_c}(\boldsymbol a))$ is the vector of the $n_c$ $p$-order polyharmonic splines defined at $n_c$ control points $\boldsymbol c^1, \ldots, \boldsymbol c^{n_c}$ and which are such that, for $i = 1, \ldots, n_c$
\begin{eqnarray}
\label{eq30}
\begin{split}
{b}_i(\boldsymbol a) = \Vert \boldsymbol a - \boldsymbol c^i \Vert ^p \quad \textrm {if $p$ is odd,} \\
{b}_i(\boldsymbol a) = \Vert \boldsymbol a - \boldsymbol c^i \Vert ^p \, \log(\Vert \boldsymbol a - \boldsymbol c^i \Vert) \quad \textrm {if $p$ is even} .
\end{split}
\end{eqnarray}
For the odd and even cases, if $\boldsymbol a =  \boldsymbol c^i$ then ${b}_i(\boldsymbol a) = 0$. In Eq.~(\ref{eq29}), $\boldsymbol{w}$ is the weight vectors in $\mathbb{R}^{n_c}$. Classically, this vector should be calculated by solving the set of $n_c$ equations
\begin{eqnarray}
\label{eq31}
\mathcal{D}(\boldsymbol{c}^i) = \boldsymbol{w}^T \boldsymbol{b}(\boldsymbol{c}^i), \quad i=1, \ldots, n_c\,,
\end{eqnarray}
yielding $n_c$ independent equations allowing vectors $\boldsymbol{w}$ to be calculated.
The function $\Psi_T ^{\rm{reg}}(\boldsymbol{a})$ in Eq.~(\ref{eq18}) is then approximated by the function $\Psi_T ^{\rm mps}(\boldsymbol{a})$ defined by
\begin{eqnarray}
\label{eq31c}
\Psi_T ^{\rm mps}(\boldsymbol{a}) =  \frac{ \boldsymbol{w}^T \boldsymbol{b}(\boldsymbol a) }{T} - \log( \mbox{1\hspace{-.25em}l}_{ \mathcal{C}_a } ^{\rm{reg}} (\boldsymbol{a})) \, ,
\end{eqnarray}
whose gradient is
\begin{eqnarray}
\label{eq31d}
\boldsymbol{\nabla}_{\boldsymbol{a}} \Psi_T ^{\rm mps}(\boldsymbol{a}) =   \frac{\boldsymbol{\nabla}_{\boldsymbol{a}} \mathcal{D}^{\rm mps}(\boldsymbol{a})}{T} + \frac{\boldsymbol{\nabla}_{\boldsymbol{a}}\mbox{1\hspace{-.25em}l}_{ \mathcal{C}_a } ^{\rm{reg}} (\boldsymbol{a}) }{\mbox{1\hspace{-.25em}l}_{ \mathcal{C}_a } ^{\rm{reg}} (\boldsymbol{a})} ,
\end{eqnarray}
with
\begin{eqnarray}
\label{eq32}
\boldsymbol{\nabla}_{\boldsymbol{a}} \mathcal{D}^{\rm mps}(\boldsymbol{a}) =  [\boldsymbol{\nabla}_{\boldsymbol{a}} \boldsymbol{b}(\boldsymbol a)] \boldsymbol{w} ,
\end{eqnarray}
where
\begin{eqnarray}
\label{eq33}
\begin{split}
[\boldsymbol{\nabla}_{\boldsymbol{a}} \boldsymbol{b}(\boldsymbol a)]_{i j} = p( a_i -  c^j_i) \Vert \boldsymbol a - \boldsymbol c^j \Vert ^{p-2} \quad \textrm {if $p$ is odd,} \\
[\boldsymbol{\nabla}_{\boldsymbol{a}} \boldsymbol{b}(\boldsymbol a)]_{i j} = ( a_i - c^j_i) (1+p\log(\Vert \boldsymbol a - \boldsymbol c^j \Vert)) \Vert \boldsymbol a - \boldsymbol c^j \Vert ^{p-2} \quad \textrm {if $p$ is even.} 
\end{split}
\end{eqnarray}
It should be noted that for $p=1$, the approximated gradient is not defined at the control points. Then only the case $p \geq 2$ will be considered.
It can be shown that as $\Vert\boldsymbol{a}\Vert \rightarrow +\infty$, function $\mathcal{D}^{\rm mps}(\boldsymbol{a})$  in Eq.(\ref{eq29}) is such that
\begin{eqnarray}
\label{eq33a}
\begin{split}
\mathcal{D}^{\rm mps}(\boldsymbol{a}) \underset{\Vert\boldsymbol{a}\Vert \rightarrow +\infty}\sim \Vert \boldsymbol a \Vert ^p \sum_{i=1}^{n_c} {w_i}  \quad \textrm {if $p$ is odd,}  \\
\mathcal{D}^{\rm mps}(\boldsymbol{a}) \underset{\Vert\boldsymbol{a}\Vert \rightarrow +\infty}\sim \Vert \boldsymbol a \Vert ^p \log(\Vert \boldsymbol a\Vert)\sum_{i=1}^{n_c} {w_i}  \quad \textrm {if $p$ is even.} \\
\end{split}
\end{eqnarray}
Then since the limit of $\log( \mbox{1\hspace{-.25em}l}_{ \mathcal{C}_a } ^{\rm{reg}} (\boldsymbol{a}))$ as $\Vert\boldsymbol{a}\Vert \rightarrow +\infty$ is $-\infty$ if the domain $\mathcal{C}_a$ is bounded and zero otherwise, a sufficient condition for Eq.(\ref{eq9a}) holds is $\sum_{i=1}^{n_c} {w_i} \geq 0$. This condition condition can be enforced by imposing the following condition 
\begin{eqnarray}
\label{eq31bc}
\sum_{i=1}^{n_c} {w_i} = \epsilon >0 \, ,
\end{eqnarray}
when solving Eq.(\ref{eq31}) (introduction of a Lagrange multiplier).
The verification of Eq.(\ref{eq9c}) is straightforward using the preceding remark and the continuity of function $\boldsymbol{b}(\boldsymbol a)$.
Concerning the condition Eq.(\ref{eq9b}), we have $\int_{\mathbb{R}^{N}} \Vert \boldsymbol{\nabla}_{\boldsymbol{u}}\Psi^{\rm mps}_T(\boldsymbol{u}) \Vert \, $ $ \exp(-\Psi^{\rm mps}_T(\boldsymbol{u}))\, \textrm{d} \boldsymbol{u}  = \int_{\mathbb{R}^{N}} { \mbox{1\hspace{-.25em}l}_{ \mathcal{C}_a } ^{\rm{reg}} (\boldsymbol{u}) \Vert \boldsymbol{\nabla}_{\boldsymbol{u}}\Psi^{\rm mps}_T(\boldsymbol{u}) \Vert \,  \exp(-\mathcal{D}^{\rm mps}(\boldsymbol{u}))\, \textrm{d} \boldsymbol{u} } \leq \int_{\mathbb{R}^{N}} {\mbox{1\hspace{-.25em}l}_{ \mathcal{C}_a } ^{\rm{reg}} (\boldsymbol{u})} $ ${\Vert  (\boldsymbol{\nabla}_{\boldsymbol{u}} \mathcal{D}^{\rm mps}(\boldsymbol{u}))/T \Vert \,  \exp(-\mathcal{D}^{\rm mps}(\boldsymbol{u})) \, \textrm{d} \boldsymbol{u} } + \int_{\mathbb{R}^{N}} {\Vert {\boldsymbol{\nabla}_{\boldsymbol{u}}\mbox{1\hspace{-.25em}l}_{ \mathcal{C}_a } ^{\rm{reg}} (\boldsymbol{u}) } \Vert \,  \exp(-\mathcal{D}^{\rm mps}(\boldsymbol{u}))\, \textrm{d} \boldsymbol{u} }$. Since (1) the function $\mbox{1\hspace{-.25em}l}_{ \mathcal{C}_a } ^{\rm{reg}}(\boldsymbol{a})$ generally decreases very rapidly out of $\mathcal{C}_a
$ and (2) for large values of the components of $\boldsymbol{a}$, (a) the exponential $\exp(-\mathcal{D}^{\rm mps}(\boldsymbol{a}))$ decreases rapidly under the condition Eq.(\ref{eq31bc}) and (b) the gradient $\boldsymbol{\nabla}_{\boldsymbol{a}} \mathcal{D}^{\rm mps}(\boldsymbol{a})$ is polynomial, then it can be deduced that the two integrals are finite and then Eq.(\ref{eq9c}) holds.

In this section, we have used polyharmonic splines for approximating the function $\mathcal{D}(\boldsymbol{a})$. Other radial basis functions such as Gaussian or Hardy \cite{Buhmann2003} types could also be used. The polyharmonic functions have been preferred since they do not require shape parameters to be estimated. It should be noted that, even if radial interpolation methods are well adapted for the interpolation of  multivariate functions, they suffer, like other surface response, of the \textit{curse of dimensionality}. Therefore, they are not adapted if the input research space dimension is very high.

\section{Approximate ISDE-based simulated annealing algorithm}
\label{sec5}

The approximate ISDE-based simulated annealing algorithm is similar to Algorithm \ref{algo2}, with the difference that function $\Psi_{T_k}^{\rm{reg}}$ is replaced by its polyharmonic splines approximation $\Psi_{T_k}^{\rm{mps}}$ introduced in the previous section. Furthermore after each iteration $k$, the polyharmonic splines approximation $\Psi_{T_k}^{\rm{mps}}$ is enriched by adding a new control point as the end value reached during the integration of the ISDE.
The algorithm is initialized by generating $n_c$ control points randomly in $\mathcal{C}_a$. The complete algorithm is summarized in Algorithm \ref{algo3}.

\begin{algorithm}[h!]
\textbf{INITIALIZATION:}\\
Generate $\boldsymbol{U}_{0}$ and $\boldsymbol{V}_{0}$\,\;\\
$\boldsymbol{A}^{1} = \boldsymbol{U}_{0}$\,\;\\
$\boldsymbol{U}^{1} = \boldsymbol{U}_{0}$\,\;\\
$\boldsymbol{V}^{1} = \boldsymbol{V}_{0}$\,\;\\
Generate $n_c$ control points $\boldsymbol c^1, \ldots, \boldsymbol c^{n_c}$\,\;\\
Calculate $\mathcal{D}(\boldsymbol c^1), \ldots, \mathcal{D}(\boldsymbol c^{n_c})$\,\;\\
Calculate $\boldsymbol{w}$ solving Eq.~(\ref{eq31})\,\;\\
\textbf{LOOP:}\\
	\For{$k=1,\ldots,(n_t-1)$}{
	Set temperature $T_k$\,\;\\
	Set damping matrix $[D_k] = [S_k][S_k]^T$\,\;\\
	Set increment $\Delta r_k$\,\;\\
	Set number of steps $M_k$\,\;\\
	Solve ISDE : \\
	\For{$\ell=1,\ldots,(M_k-1)$}{
		$\boldsymbol{V}^{\ell+1}   =   ([I]-\frac{\Delta r_k}{2} [D_k]) \boldsymbol{V}^{\ell}- \Delta r_k \boldsymbol{\nabla}_{\boldsymbol{u}} \Psi_{T_k}^{\rm{mps}} (\boldsymbol{U}^{\ell})  +[S_k] \Delta \boldsymbol{W}^{\ell+1}$\,\;\\
		$\boldsymbol{U}^{\ell+1}   =   \boldsymbol{U}^{\ell} + \Delta r_k\, \boldsymbol{V}^{\ell+1}$\,\;\\
		}
		$\boldsymbol{A}^{k+1} = \boldsymbol{U}_{M_k}$\,\;\\
		$\boldsymbol{U}^{1} \leftarrow \boldsymbol{U}_{M_k}$\,\;\\
		$\boldsymbol{V}^{1} \leftarrow \boldsymbol{V}_{M_k}$\,\;\\
		$n_c \leftarrow n_c + 1$\,\;\\
		Add $\boldsymbol c^{n_c} = \boldsymbol{A}^{k+1}$ to the list of control points\,\;\\
		Calculate $\mathcal{D}(\boldsymbol c^{n_c})$\,\;\\
		Update $\boldsymbol{w}$ solving Eq.~(\ref{eq31})\,\;\\
}
\caption{Approximate ISDE-based simulated annealing algorithm}
\label{algo3}
\end{algorithm}
This procedure allows the added control points to be distributed according to the temperature-dependent distribution. This algorithm can be parallelized by integrating simultaneously several ISDEs with different randomly generated initial conditions at each iteration $k$.
\section{Application}
\label{sec6}
\subsection{Ackley's test function}
The objective here is to illustrate and show the efficiency of both Algorithms \ref{algo2} and \ref{algo3}.  For this purpose, we use the  Ackley's function which is a classical function for testing multivariate optimization algorithms. This function is defined by
\begin{eqnarray}
\label{eq36}
 \mathcal{D}(\boldsymbol{a}) = -20\exp(-\frac{0.2}{\sqrt N} \Vert \boldsymbol{a} \Vert )- \exp(\frac{1}{N}\sum_{i=1}^N(cos(2\pi {a_i})) + \exp(1) +20,
\end{eqnarray}
with the following constraints
\begin{eqnarray}
\label{eq37}
 -5<a_i<5,\quad  i=1,\ldots,N\,.
\end{eqnarray}
For this case the global maximum is reached at
\begin{eqnarray}
\label{eq38}
 a_i \simeq 0,\quad  i=1,\ldots,N\,.
\end{eqnarray}
For $N=1$ and $N=2$ the cost functions are plotted on Figs.~\ref{fig1} and~\ref{fig2}.
\begin{figure}[h!]
\begin{center}
\includegraphics[width=6.5cm]{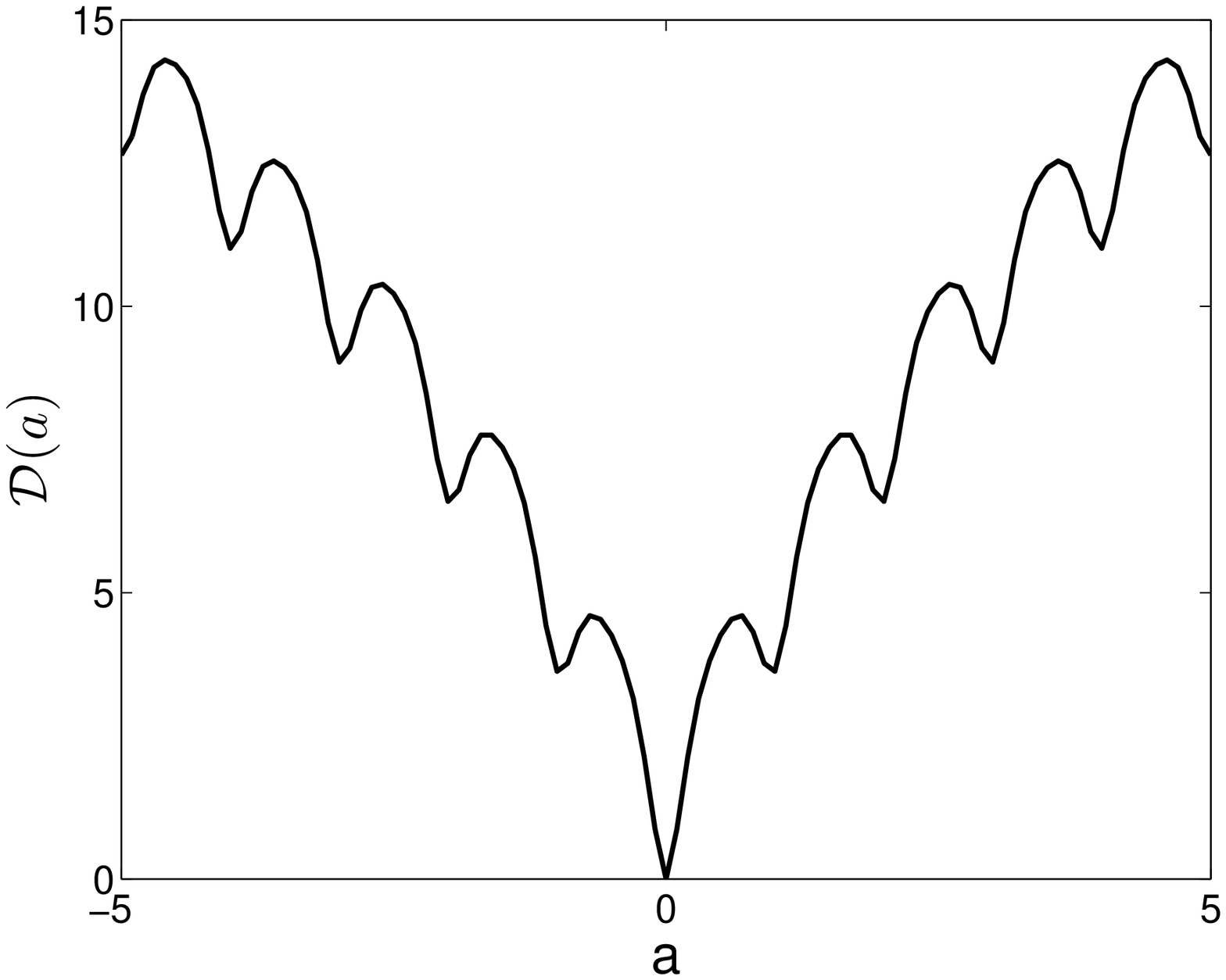}
\caption{Function $\boldsymbol{a} \mapsto \mathcal{D}(\boldsymbol{a})$ for $N=1$. \label{fig1}}
\end{center}
\end{figure}
\begin{figure}[h!]
\begin{center}
\includegraphics[width=6.5cm]{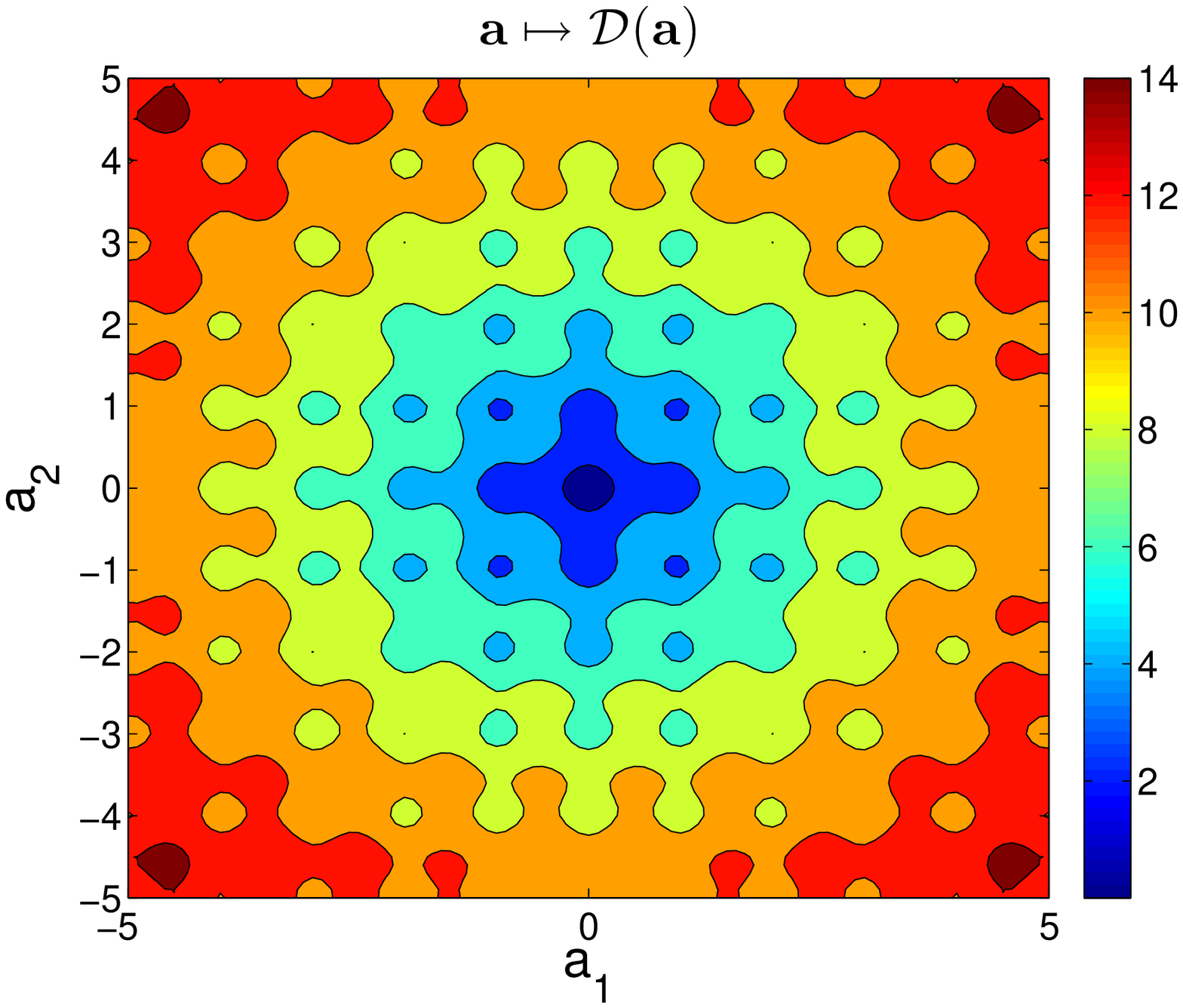}
\caption{Function $\boldsymbol{a} \mapsto \mathcal{D}(\boldsymbol{a})$ for $N=2$. \label{fig2}}
\end{center}
\end{figure}
The number of local minima increases rapidly with dimension.
\subsubsection{Exact ISDE-based simulated annealing algorithm}

For case $N=2$, the ISDE-based simulated annealing sequence, described in Algorithm \ref{algo2}, is generated. The total number of ISDE integrations is $n_t = 500$. For each ISDE, there are $M_k=40$ steps. The temperature decreasing law is chosen as
\begin{eqnarray}
\label{eq39}
 T_k = T_1\,\,exp(-\beta\, k) + b.
\end{eqnarray}
in which $T_1 = 36.7$, $\beta = 2.0\times 10^{-2}$, $b = 3.51\times 10^{-2}$. 
This function is represented in Fig.~\ref{fig3}.
\begin{figure}[h!]
\begin{center}
\includegraphics[width=6.5cm]{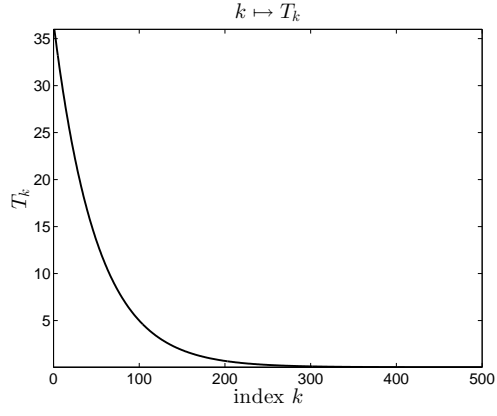}
\caption{Function $k \mapsto T_k$. \label{fig3}}
\end{center}
\end{figure}
To take into account the constraints, the regularized indicator function in Eq.(\ref{eq17}) with $\alpha_1=\alpha_2 =0.3$. Figure~\ref{fig4} shows the values $\boldsymbol{A}_1, \ldots, \boldsymbol{A}_{n_t}$ generated using Algorithm \ref{algo2}.
\begin{figure}[h!]
\begin{center}
\includegraphics[width=6.5cm]{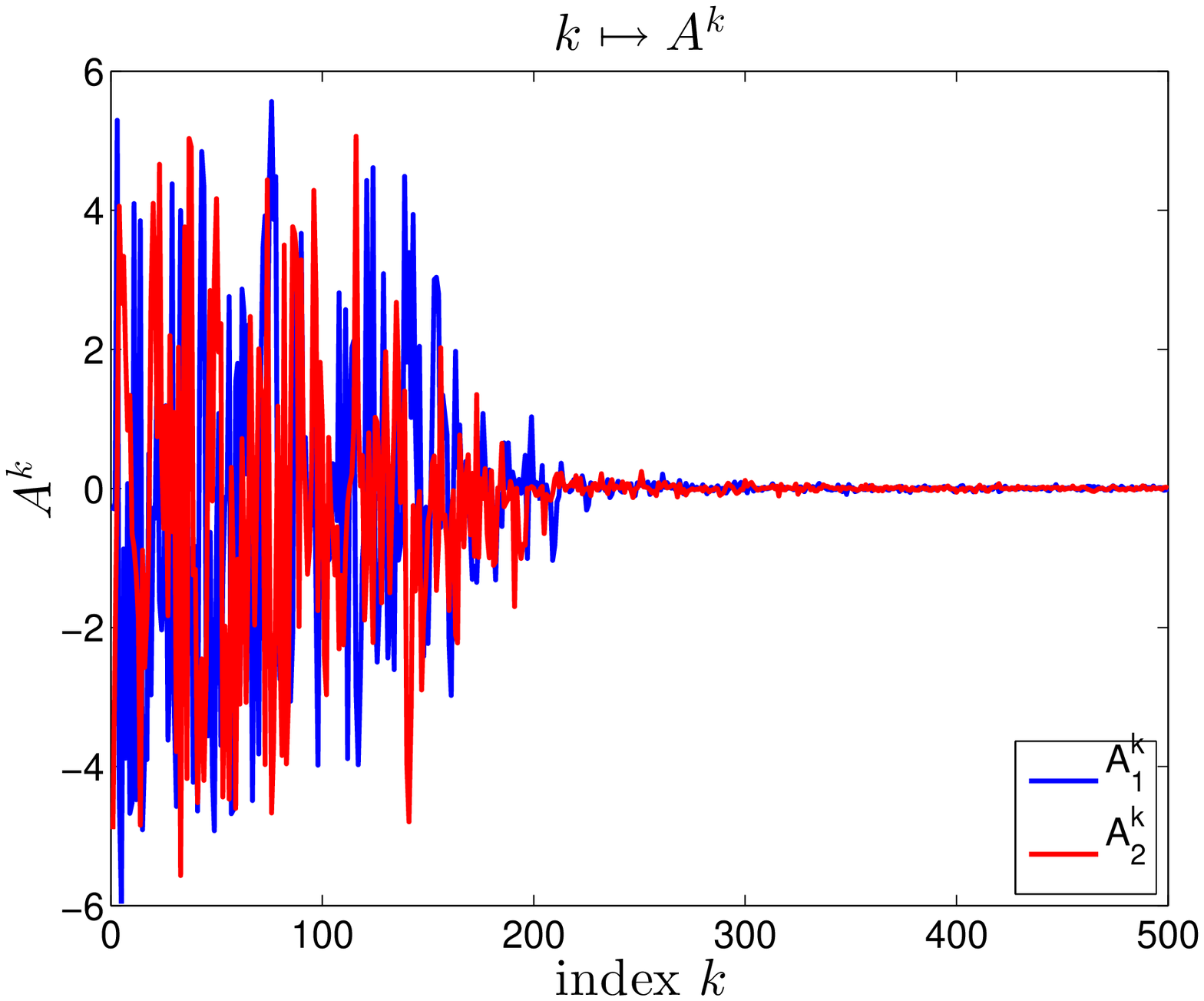}
\includegraphics[width=6.5cm]{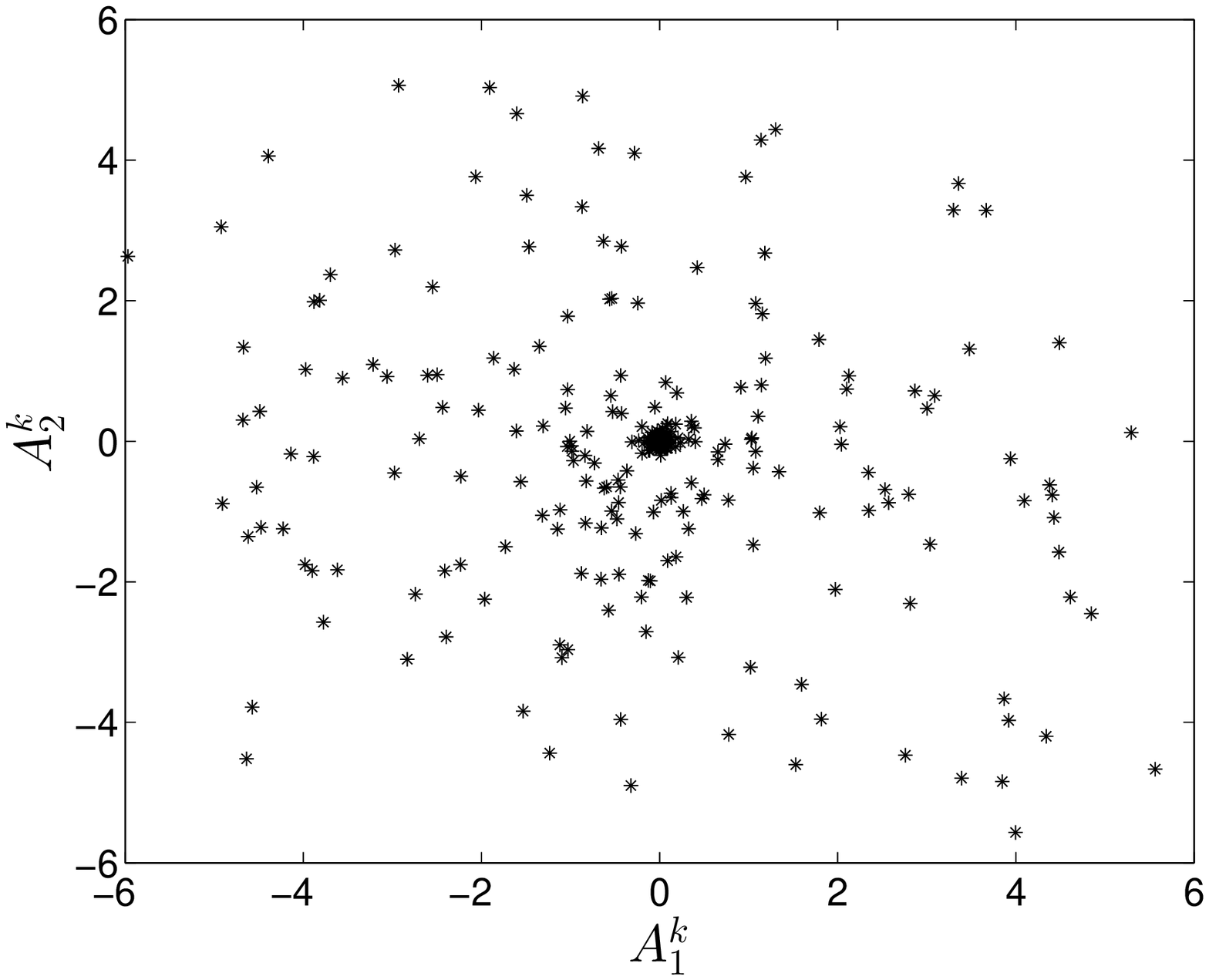}
\caption{ISDE-based algorithm, $N=2$, functions $k \mapsto A_1^k$ and $k \mapsto A_2^k$. \label{fig4}}
\end{center}
\end{figure}
It can be seen the capacity of the algorithm to escape from local minima and to converge to the global minimum. There are still fluctuations at the end of the simulation since the end-temperature is not zero. 

Let's now increase the dimension to $N=256$. The parameter of the temperature decreasing law are $T_1 = 2.5$, $\beta = 2.0\times 10^{-2}$, $b = 5.1\times 10^{-3}$.
Figure~\ref{fig4b} shows the values $\boldsymbol{A}_1, \ldots, \boldsymbol{A}_{n_t}$ generated using Algorithm \ref{algo2}).
\begin{figure}[h!]
\begin{center}
\includegraphics[width=6.5cm]{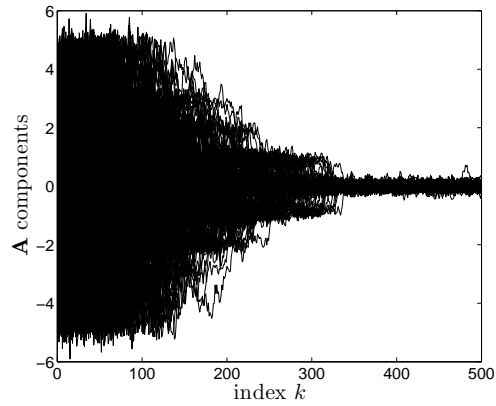}
\caption{ISDE-based algorithm, $N=200$, functions $k \mapsto A_i^k$ for $i=1,\ldots,200$. \label{fig4b}}
\end{center}
\end{figure}
It can be seen in this figure that all the components converge to the global minimum. To show the efficiency of this algorithm, these results are compared with the ones obtained using the classical simulated annealing algorithm  (\ref{algo1}). The same temperature law, and the same increment size calculation method (based on the Hessian of cost function) has been used. Nevertheless, since the calculations of the gradient in Algorithm (\ref{algo2}) is more expensive than the simple function evaluations in Algorithm (\ref{algo1}), $M_k=170$ steps for each MH sequence has been set in order to have the same calculation time. The results of the optimisation using the classical simulated annealing algorithm are shown in Fig.~\ref{fig4c}.
\begin{figure}[h!]
\begin{center}
\includegraphics[width=6.5cm]{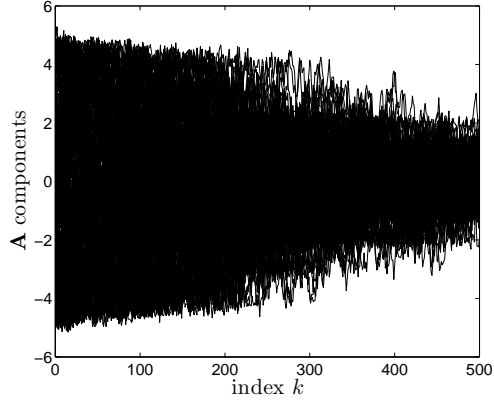}
\caption{Classical SA algorithm,  $N=200$, functions $k \mapsto A_i^k$ for $i=1,\ldots,200$. \label{fig4c}}
\end{center}
\end{figure}
It can be seen in this figure that the result obtained using the ISDE based algorithm are better than those obtained using the classical SA algorithm. 
\subsubsection{Approximate ISDE-based simulated annealing algorithm}
This section illustrates Algorithm \ref{algo3} for which the cost function is approximated using a polyharmonic spline approximation at order $p=2$. The initial number of control points is $n_c = 140$. The number of enrichments is $n_t = 500$. For each ISDE, there are $M_k=40$ steps. It should be noted that for this simple application the computation time is larger using the adaptative Algorithm \ref{algo3} than using the Algorithm \ref{algo2} (without approximation of the cost function). The objective here is just to validate the methodology. The computational gain of Algorithm \ref{algo3} will be illustrated in the next application which consists in a Finite Element simulation.
For case $N=2$, Fig.~\ref{fig5} shows the values $\boldsymbol{A}_1, \ldots, \boldsymbol{A}_{n_t}$. It can be seen again that the algorithm converges rapidly to the global minimum. 
\begin{figure}[h!]
\begin{center}
\includegraphics[width=6.5cm]{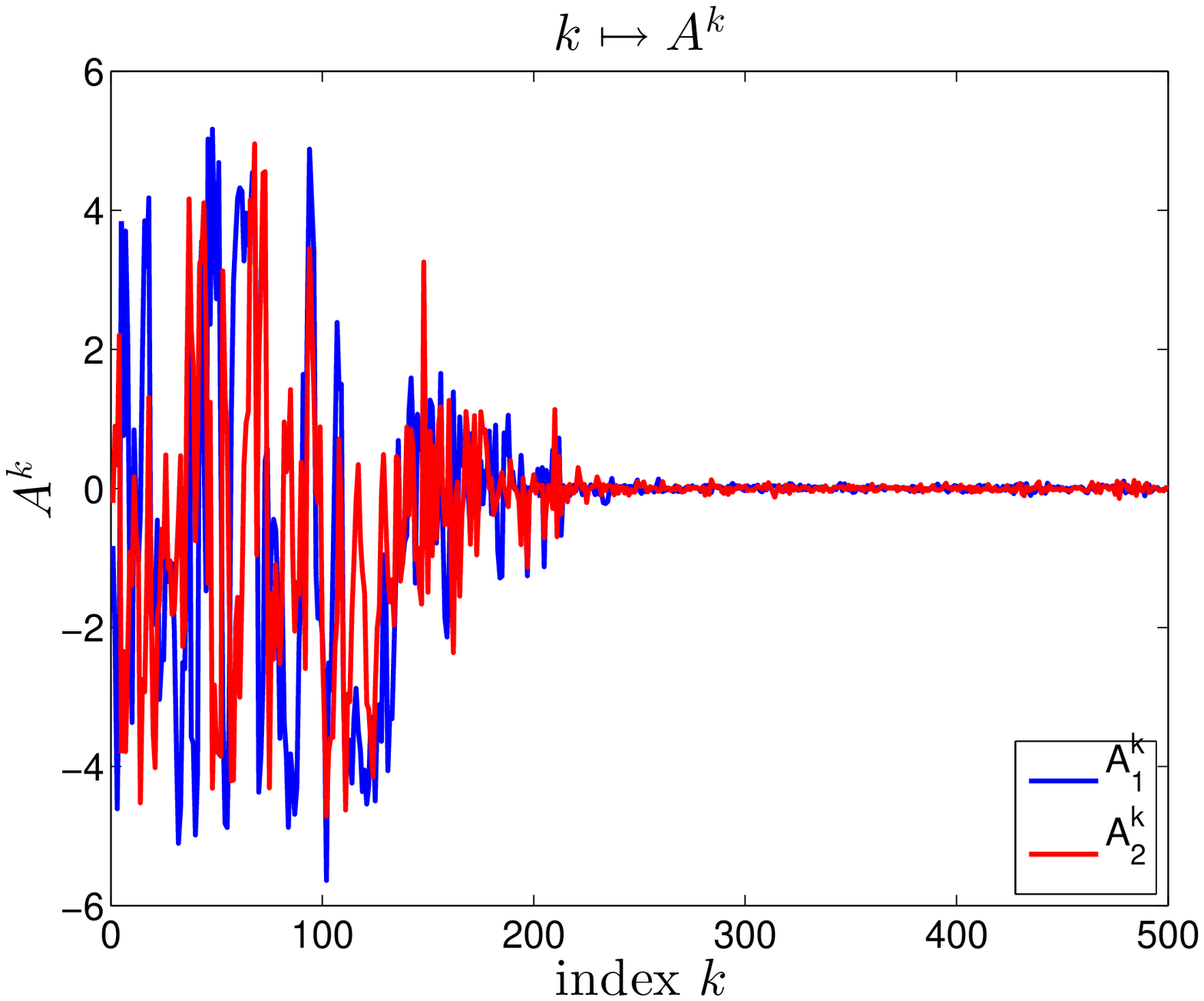}
\includegraphics[width=6.5cm]{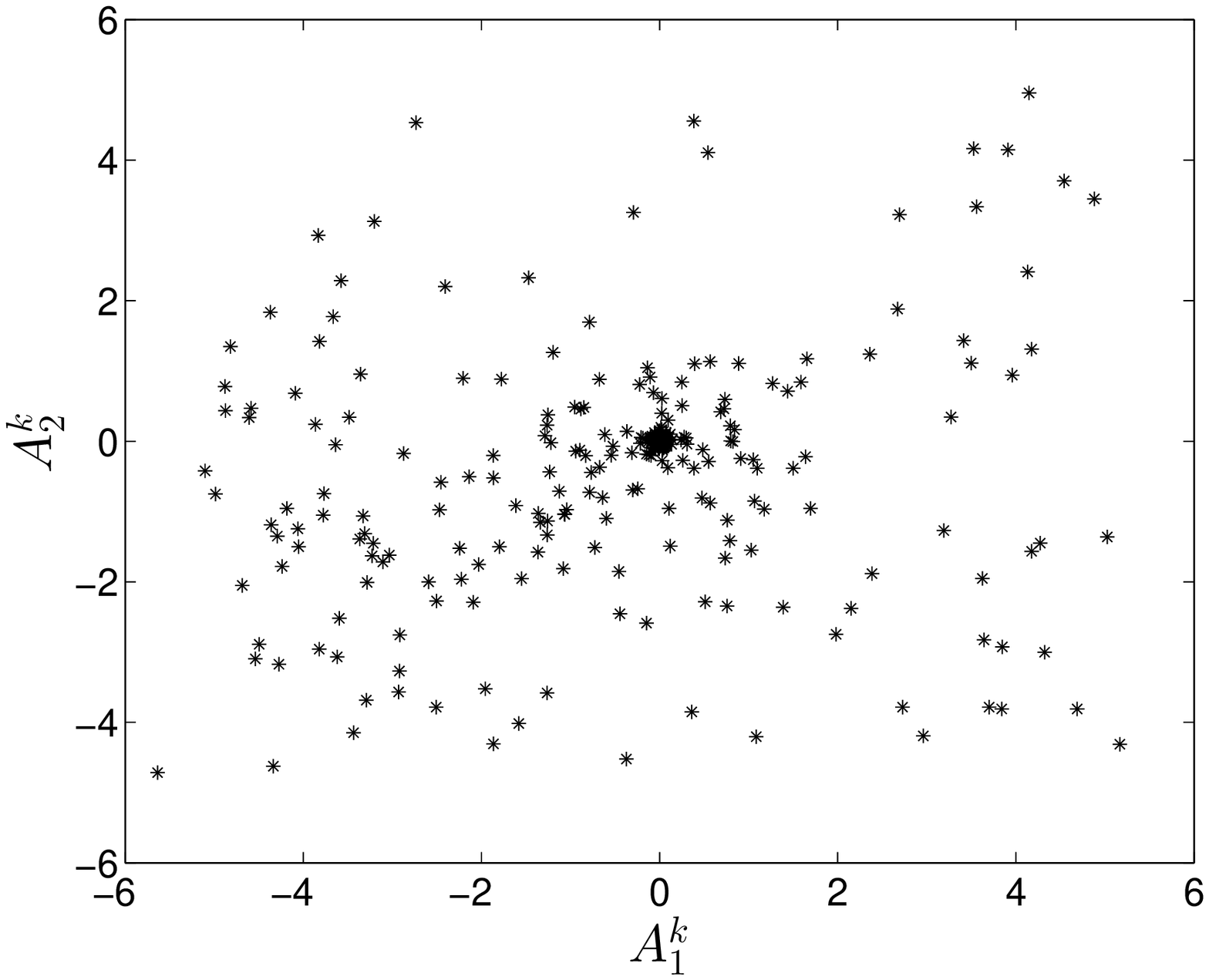}
\caption{Approximate ISDE-based algorithm, $N=2$, functions $k \mapsto A_1^k$ and $k \mapsto A_2^k$. \label{fig5}}
\end{center}
\end{figure}
The approximate cost function obtained after the $500$ enrichments is plotted on Fig.~\ref{fig6}. Compared with Fig.~\ref{fig2}, it can be seen a good approximation of the cost function is the visited regions.
\begin{figure}[h!]
\begin{center}
\includegraphics[width=6.5cm]{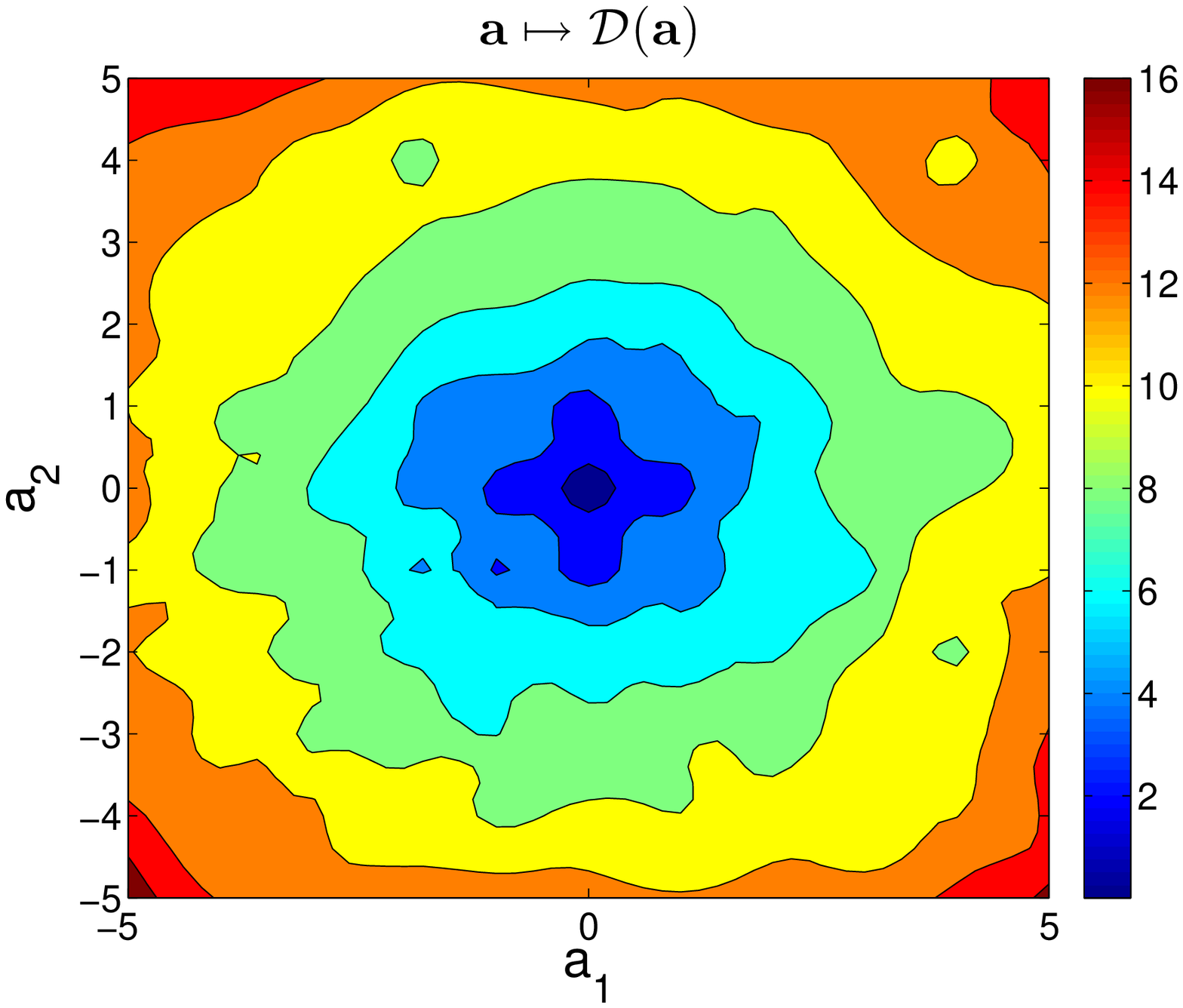}
\caption{Approximate cost function $\boldsymbol{a} \mapsto \mathcal{D}^{\rm mps}(\boldsymbol{a})$ for $N=2$. \label{fig6}}
\end{center}
\end{figure}
For case $N=32$, Fig.~\ref{fig7} shows the values $\boldsymbol{A}_1, \ldots, \boldsymbol{A}_{n_t}$. It can be seen the ability of the algorithm to converge to the global minimum when increasing the dimension. 
\begin{figure}[h!]
\begin{center}
\includegraphics[width=6.5cm]{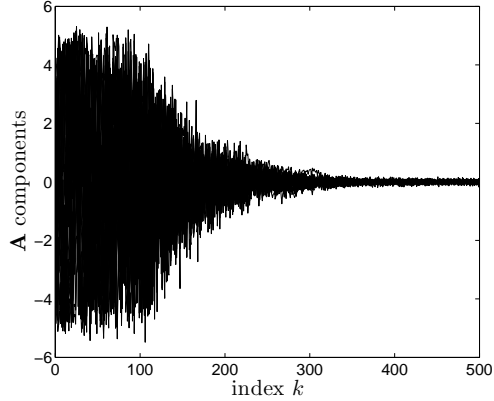}
\caption{Approximate ISDE-based algorithm, $N=32$, functions $k \mapsto A_i^k$ for $i=1,\ldots,32$. \label{fig7}}
\end{center}
\end{figure}
\subsection{Finite Element model}
We are interested in the maximal response of the 2-storey structure represented in Fig.~\ref{fig7} submitted to the soil acceleration plotted in Fig.~\ref{fig8}.
\begin{figure}[h!]
\begin{center}
\includegraphics[width=10.5cm]{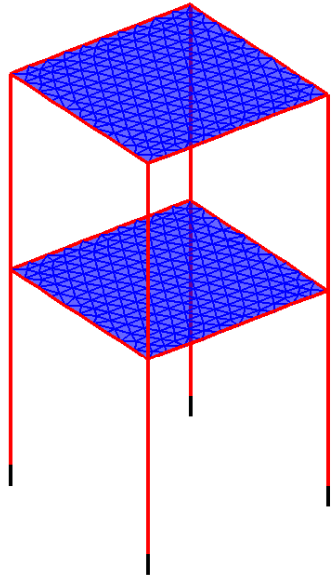}
\caption{Finite element mesh. \label{fig7}}
\end{center}
\end{figure}
This structure is made up with two $3.0\times 3.0$~m$^2$ plates (blue color), eight vertical beams (red color), eight horizontal beams (red color). The structure is linked to the ground by four linear springs.
The bottom plate has thickness $5.0\times 10^{-3}$~m, Young's modulus $6.31\times 10^{10}$~Pa, mass density $1800$~kg/m$^3$, Poisson ratio $0.29$.
The top plate has thickness $5.0\times 10^{-3}$~m, Young's modulus $6.47\times 10^{10}$~Pa, mass density $1800$~kg/m$^3$ and Poisson ratio $0.29$.
All the vertical beams have length $3$~m, diameter $3.5\times 10^{-2}$~m (circular section) Young's modulus $1.3\times 10^{11}$~Pa, mass density $7800$~kg/m$^3$ and Poisson ratio $0.3$.
All the horizontal beams have length $3$~m, diameter $3.5\times 10^{-2}$~m (circular section) Young's modulus $1.3\times 10^{11}$~Pa and Poisson ratio $0.3$.
The north, east, south and west horizontal beams have mass density $5800$~kg/m$^3$, $6800$~kg/m$^3$, $7100$~kg/m$^3$ and $6800$~kg/m$^3$ respectively.
The initial value for the springs' stiffnesses (same value for the three directions) are $k_1 = k_2 = k_3 = k_4 = 70000$~N/m. The structure is subjected to the seismic ground acceleration plotted in Fig.~\ref{fig8} along $x$-direction.
\begin{figure}[h!]
\begin{center}
\includegraphics[width=8.5cm]{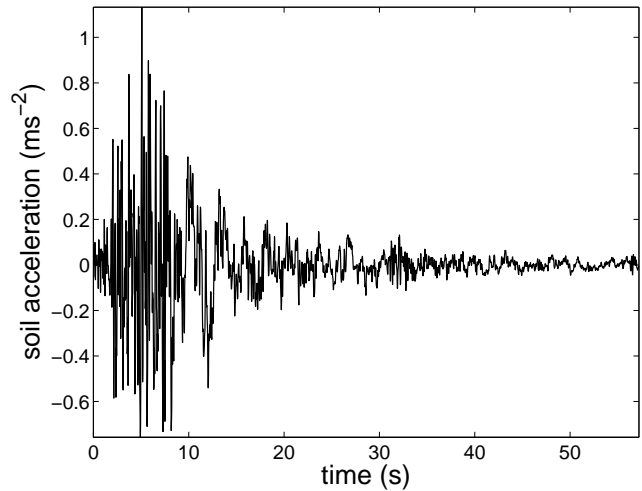}
\caption{Soil acceleration. \label{fig8}}
\end{center}
\end{figure}
We are interested in the displacement at the observation point located at $(2,2,2)$~m on the top plate. The norm of this displacement is plotted in Fig.~\ref{fig9}.
\begin{figure}[h!]
\begin{center}
\includegraphics[width=8.5cm]{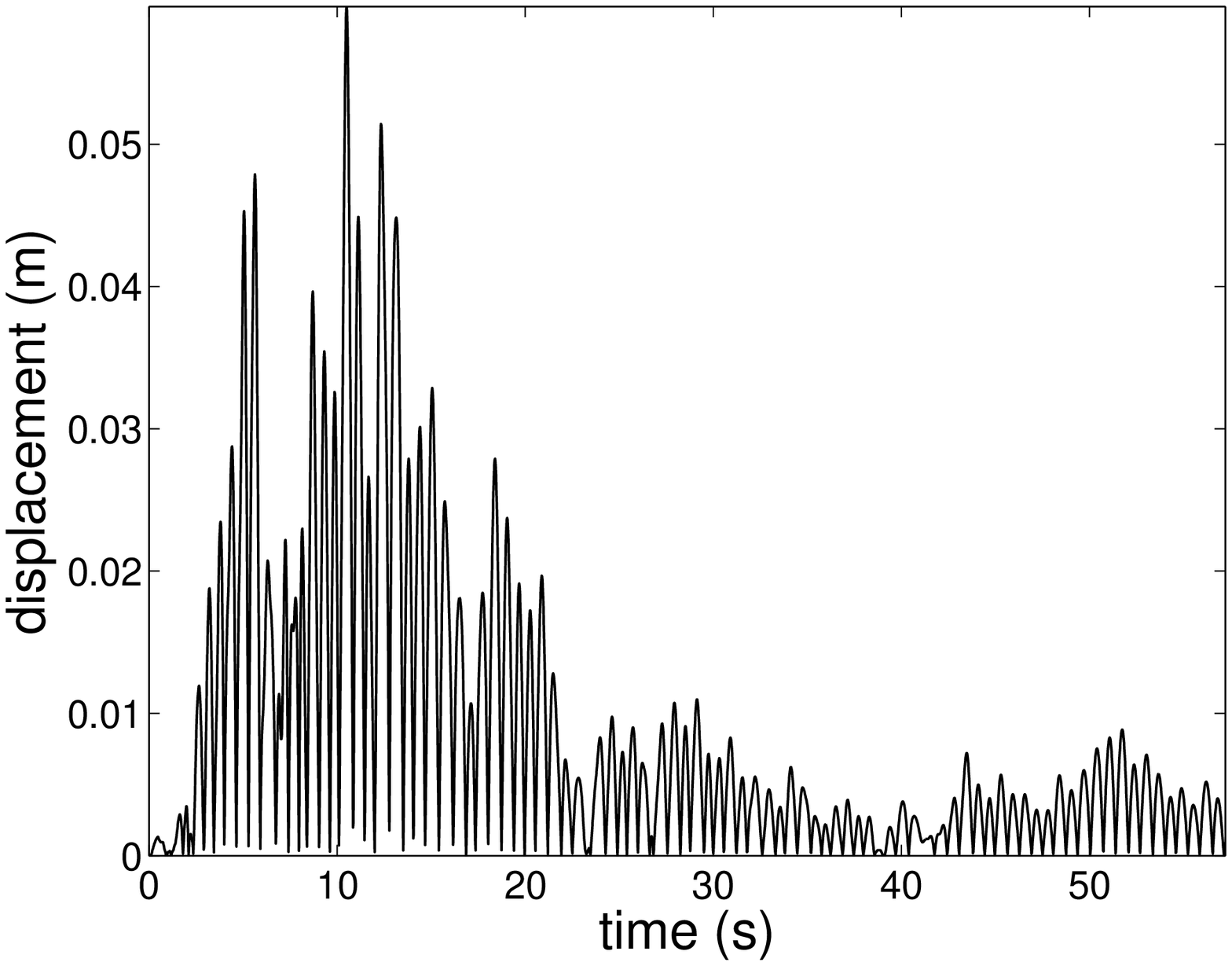}
\caption{Norm of the displacement at observation point. \label{fig9}}
\end{center}
\end{figure}
The objective is to design the stiffnesses of the springs in order to minimize the maximum value of the norm of the displacement. Then $\textbf{a} = (k_1,k_2,k_3,k_4)$.To perform this task, the approximated ISDE-based simulated algorithm \ref{algo3} is used with a polyharmonic approximation at order $p=2$. The admissible domains for the stiffnesses are $[5000, 120000]$~N/m. The initial number of control points is $n_c = 140$. The number of enrichments is $n_t = 200$. The number of steps for each ISDE between two enrichments is $M_k = 40$. The temperature decreasing law is plotted in Fig.~\ref{fig9b}. 
\begin{figure}[h!]
\begin{center}
\includegraphics[width=6.5cm]{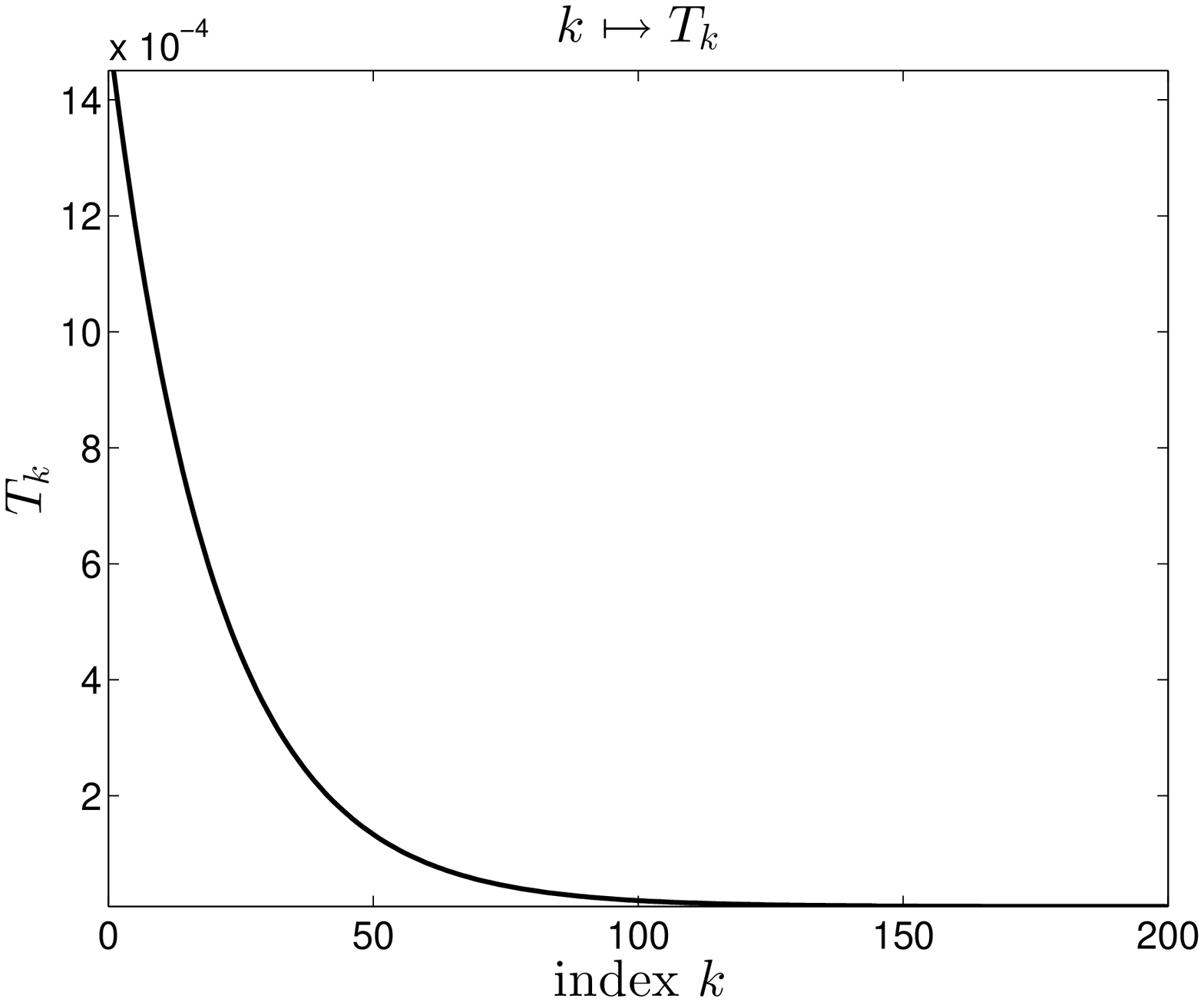}
\caption{Function $k \mapsto T_k$. \label{fig9b}}
\end{center}
\end{figure}
The optimal control point is $\textbf{a}^{\rm opt} = (19787, 18298, 59092, 70160)$~N/m. The evolution of the stiffnesses and the exact cost function during the simulation are shown in Figs.~\ref{fig10} and \ref{fig10b}. 
\begin{figure}[h!]
\begin{center}
\includegraphics[width=8.5cm]{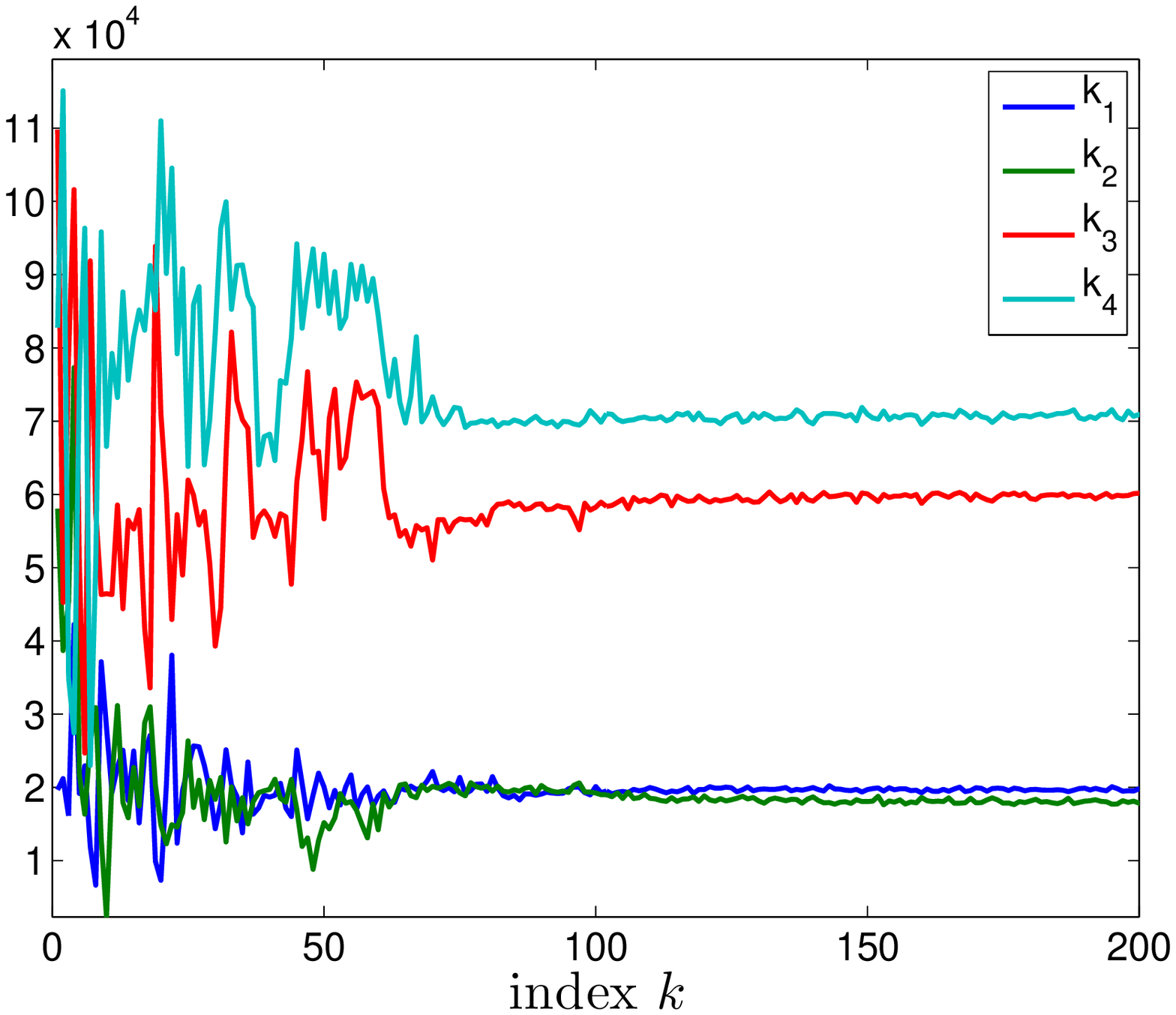}
\caption{Evolution of the stiffnesses during the simulation. \label{fig10}}
\end{center}
\end{figure}
\begin{figure}[h!]
\begin{center}
\includegraphics[width=8.5cm]{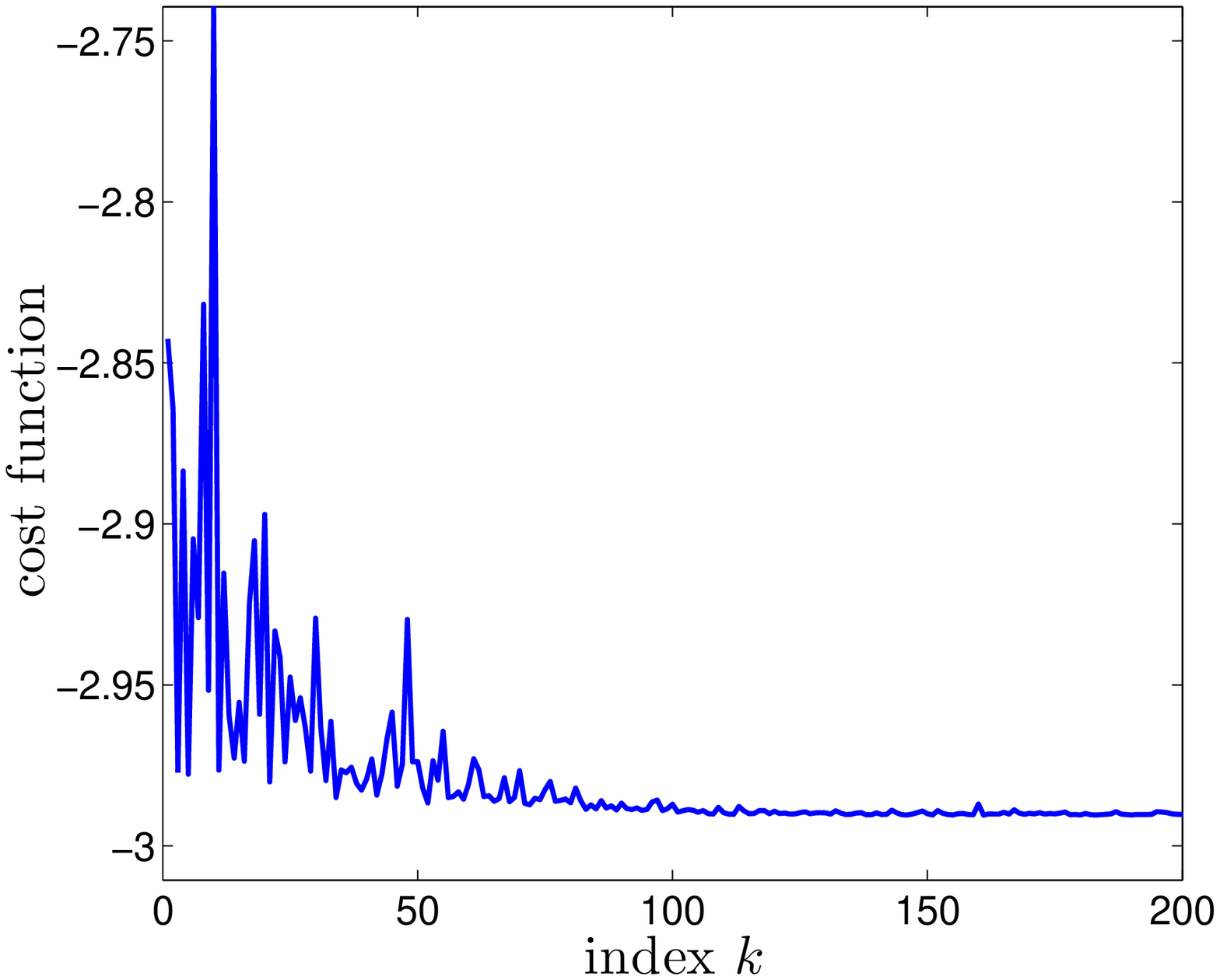}
\caption{Evolution of the exact logarithm of cost function during the simulation. \label{fig10b}}
\end{center}
\end{figure}
It can be seen it these figures that the optimum is rapidly reach using a few number of function calls.
The exact and approximate univariate variations of the cost function around the optimal value are illustrated in Fig.~\ref{fig11}.The exact and approximate bivariate variations of the cost function around the optimal value are illustrated in Figs.~\ref{fig12} and \ref{fig12b}. It can be seen in these figures that the global optimum has been found despite the presence of local minima. Furthermore, it can be seen that the polyharmonic splines provide a good approximation of the cost function around the optimal value. 
\begin{figure}[h!]
\begin{center}
\includegraphics[width=5.5cm]{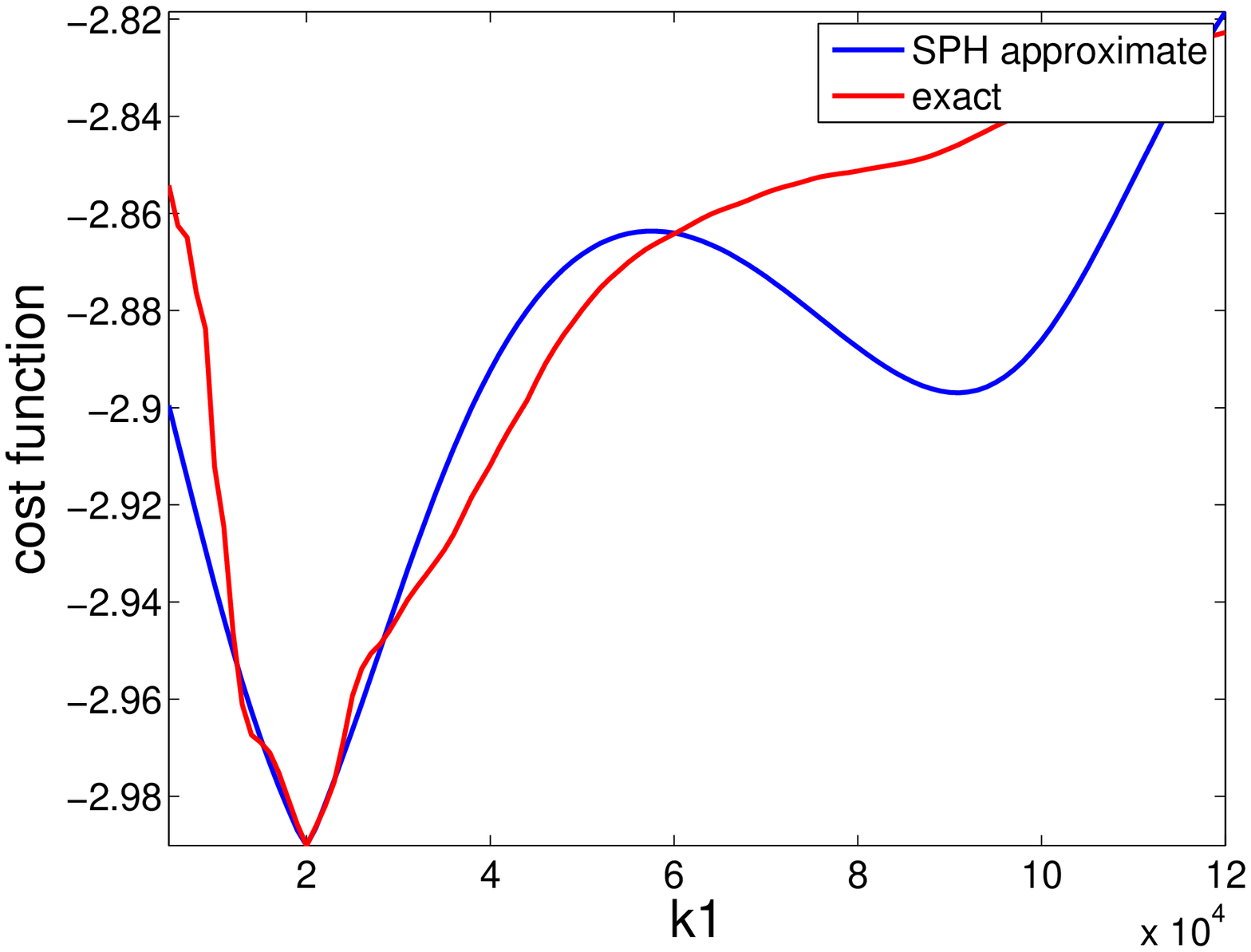}
\includegraphics[width=5.5cm]{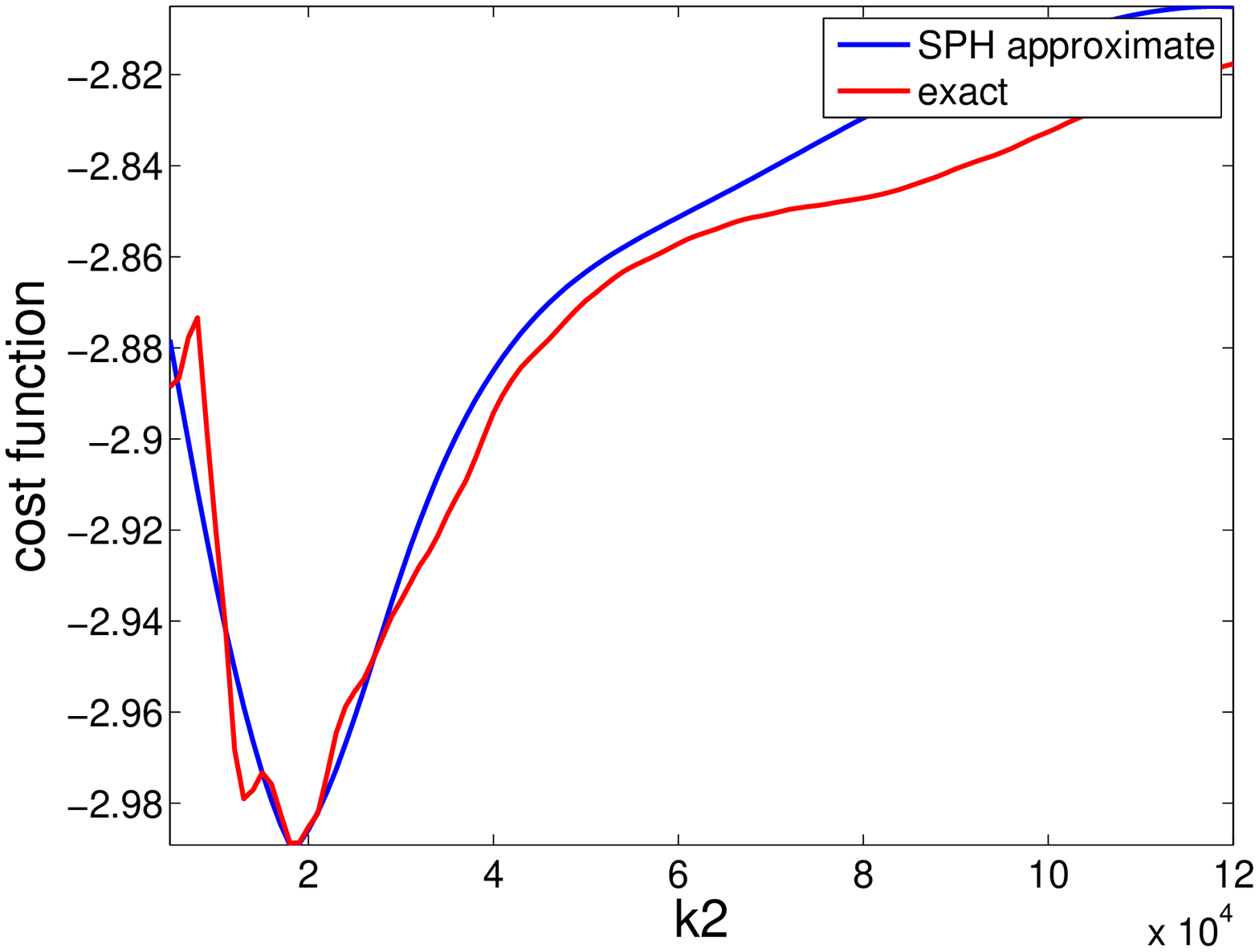}
\includegraphics[width=5.5cm]{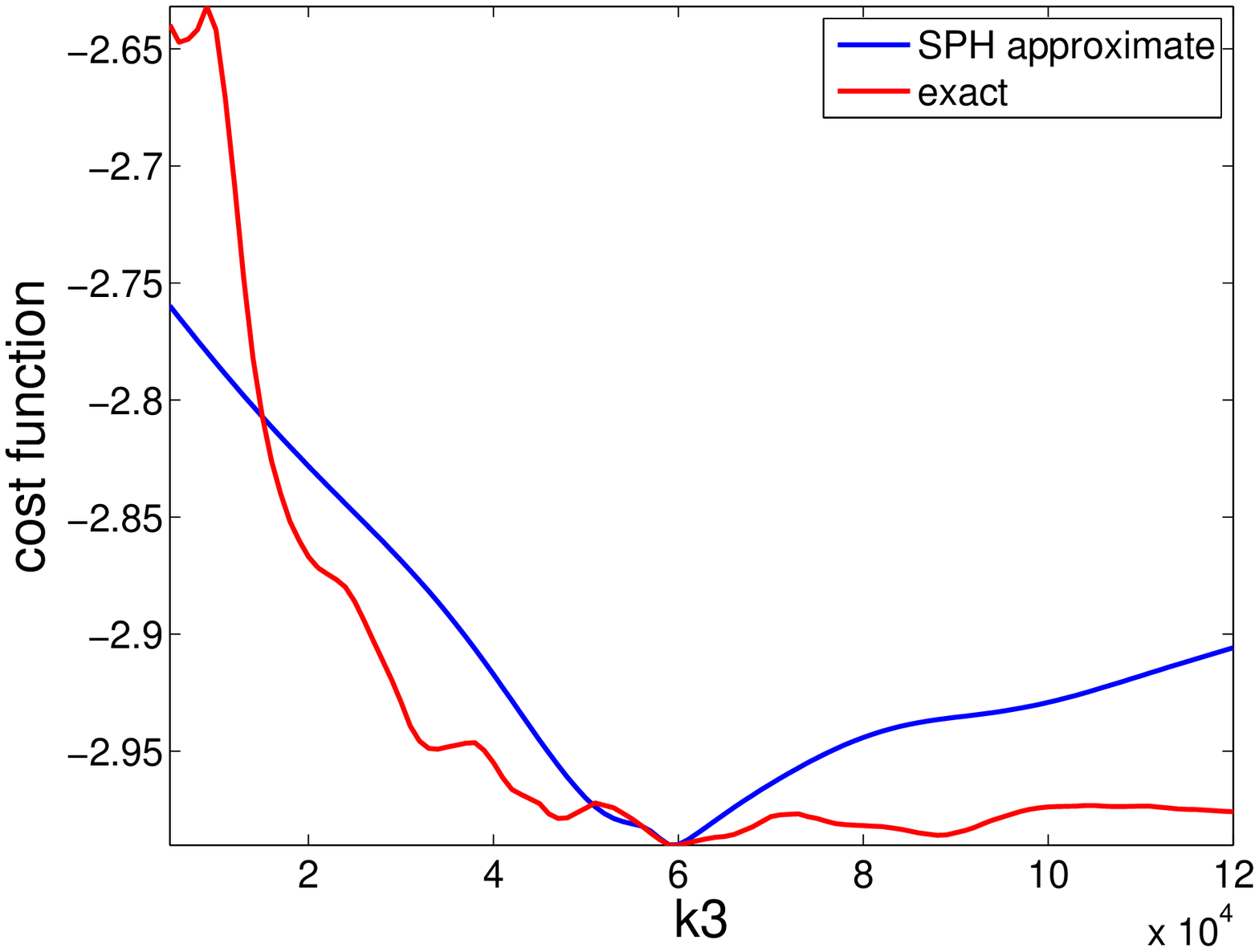}
\includegraphics[width=5.5cm]{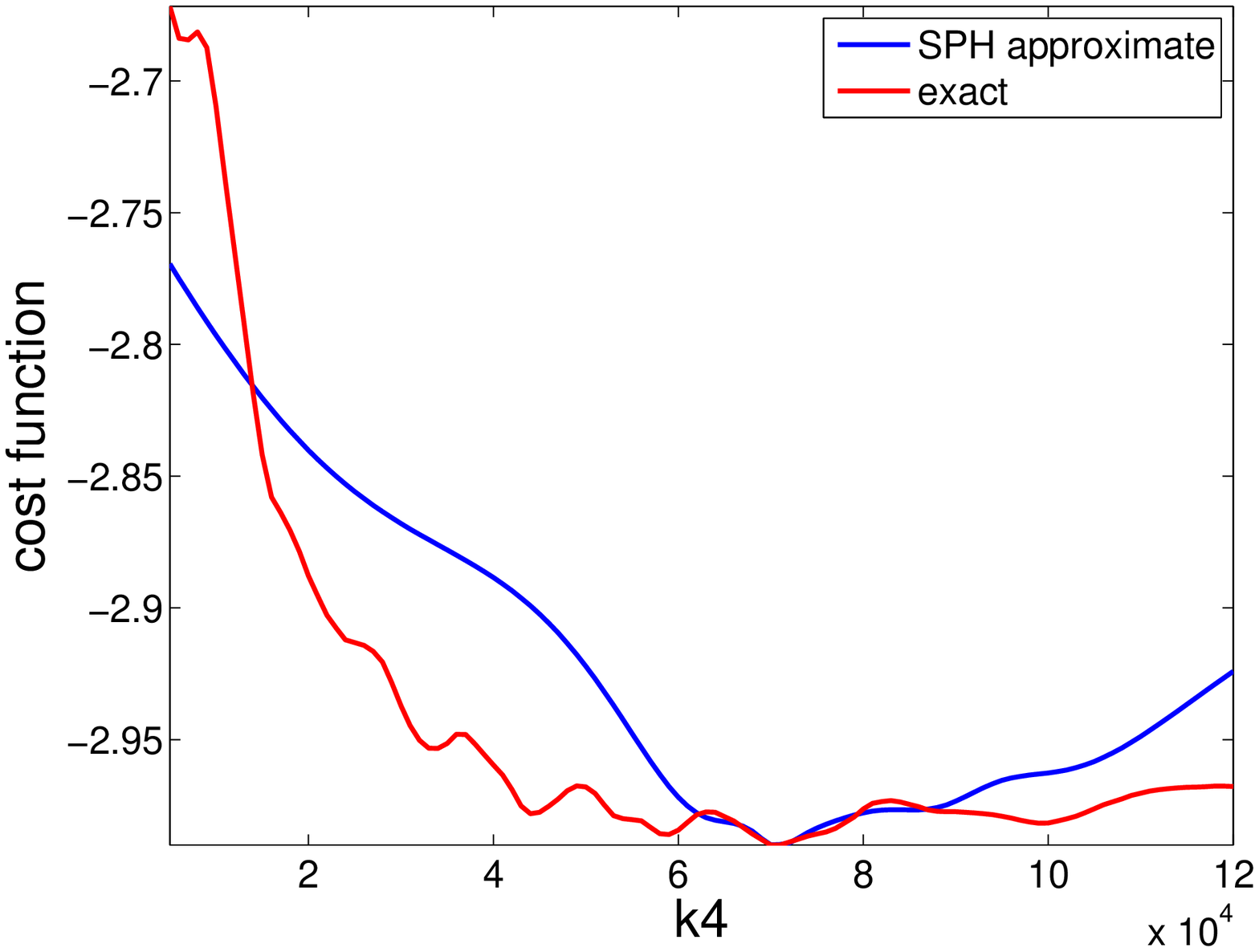}
\caption{Univariate variations of the exact and approximate logarithm of cost function around the optimal value. \label{fig11}}
\end{center}
\end{figure}
It should be noted that these figures seem to show that the algorithm has reach the global minimum but there is no guaranty unless all the research space is rigorously explored.
\begin{figure}[h!]
\begin{center}
\includegraphics[width=5.5cm]{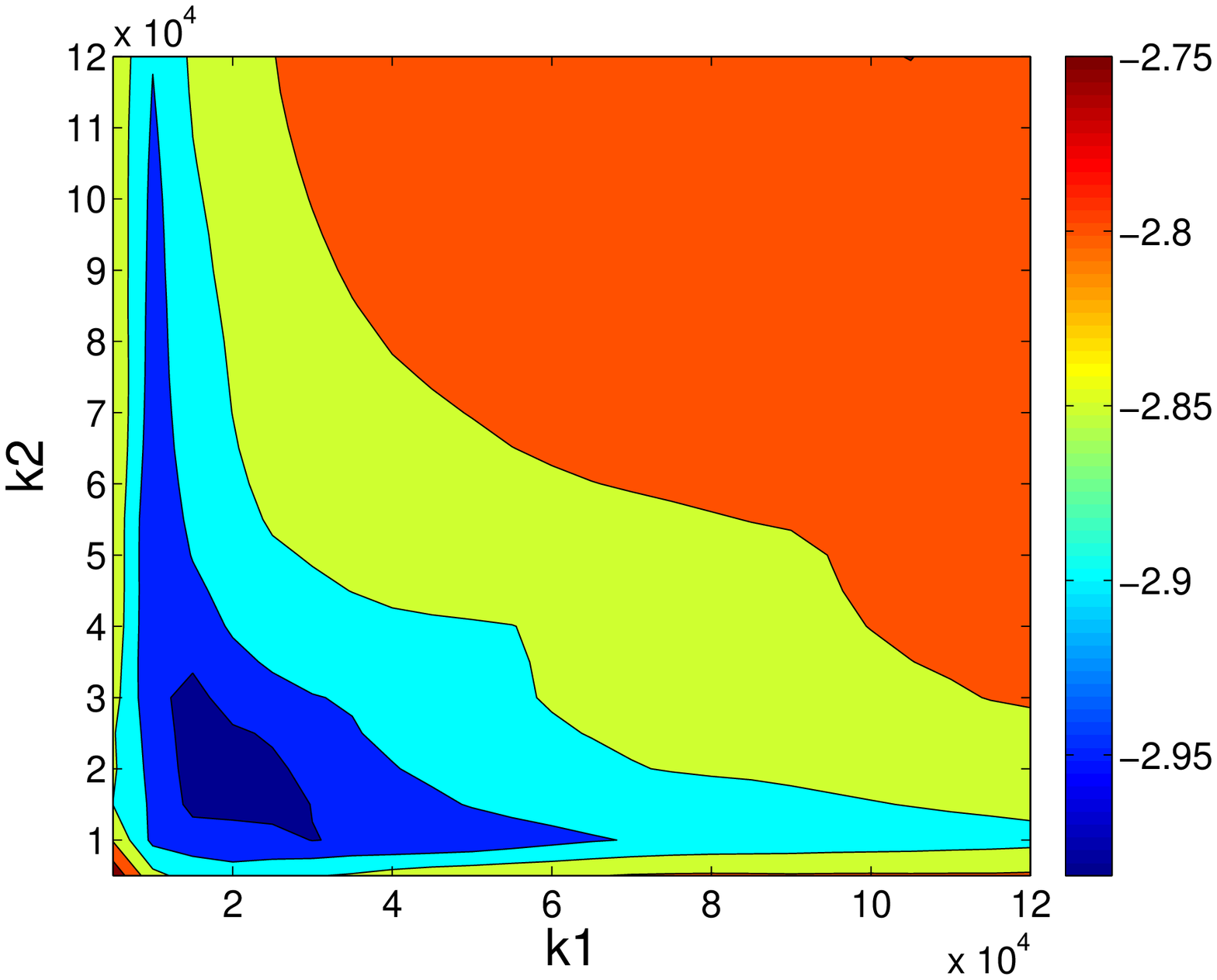}
\includegraphics[width=5.5cm]{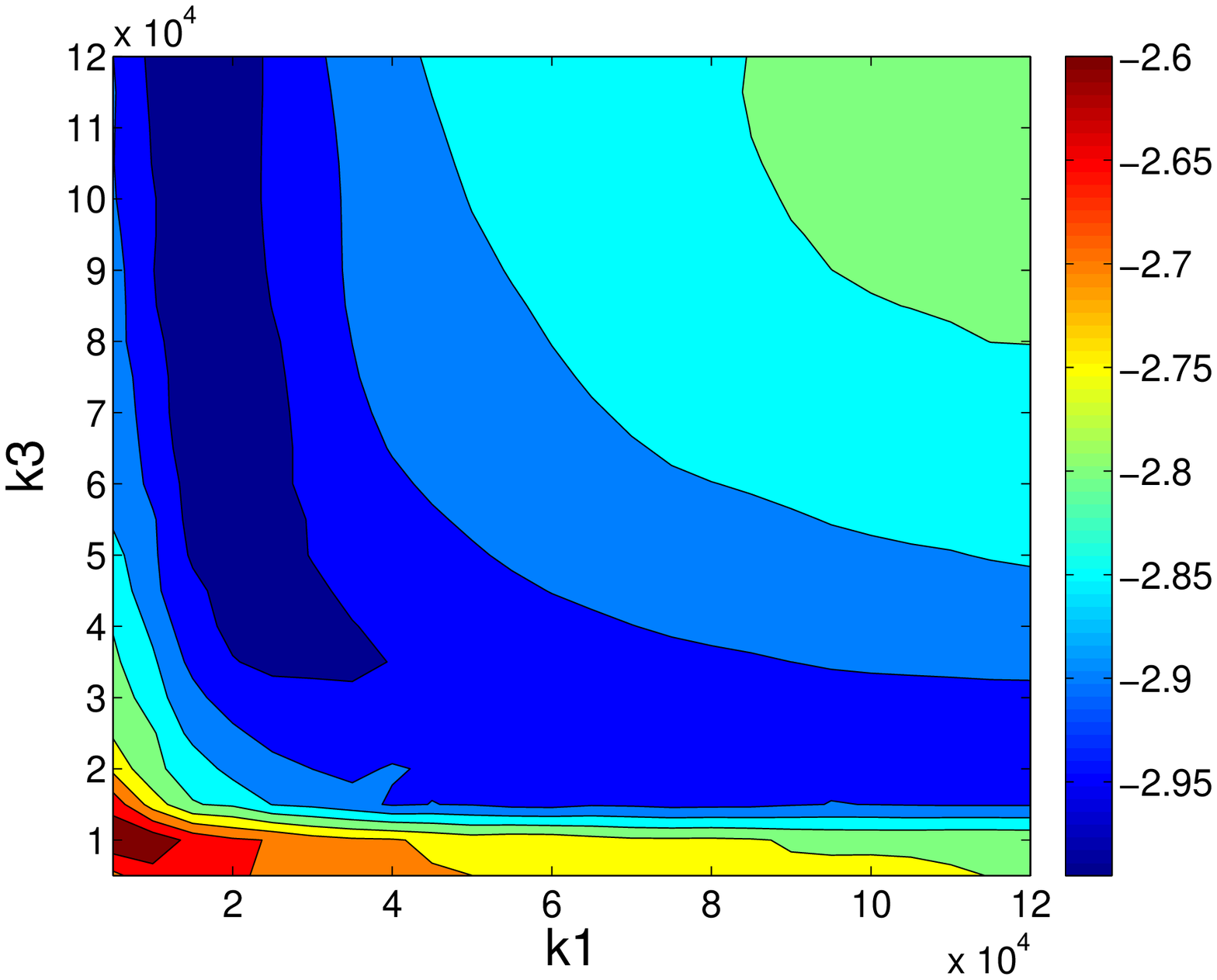}
\includegraphics[width=5.5cm]{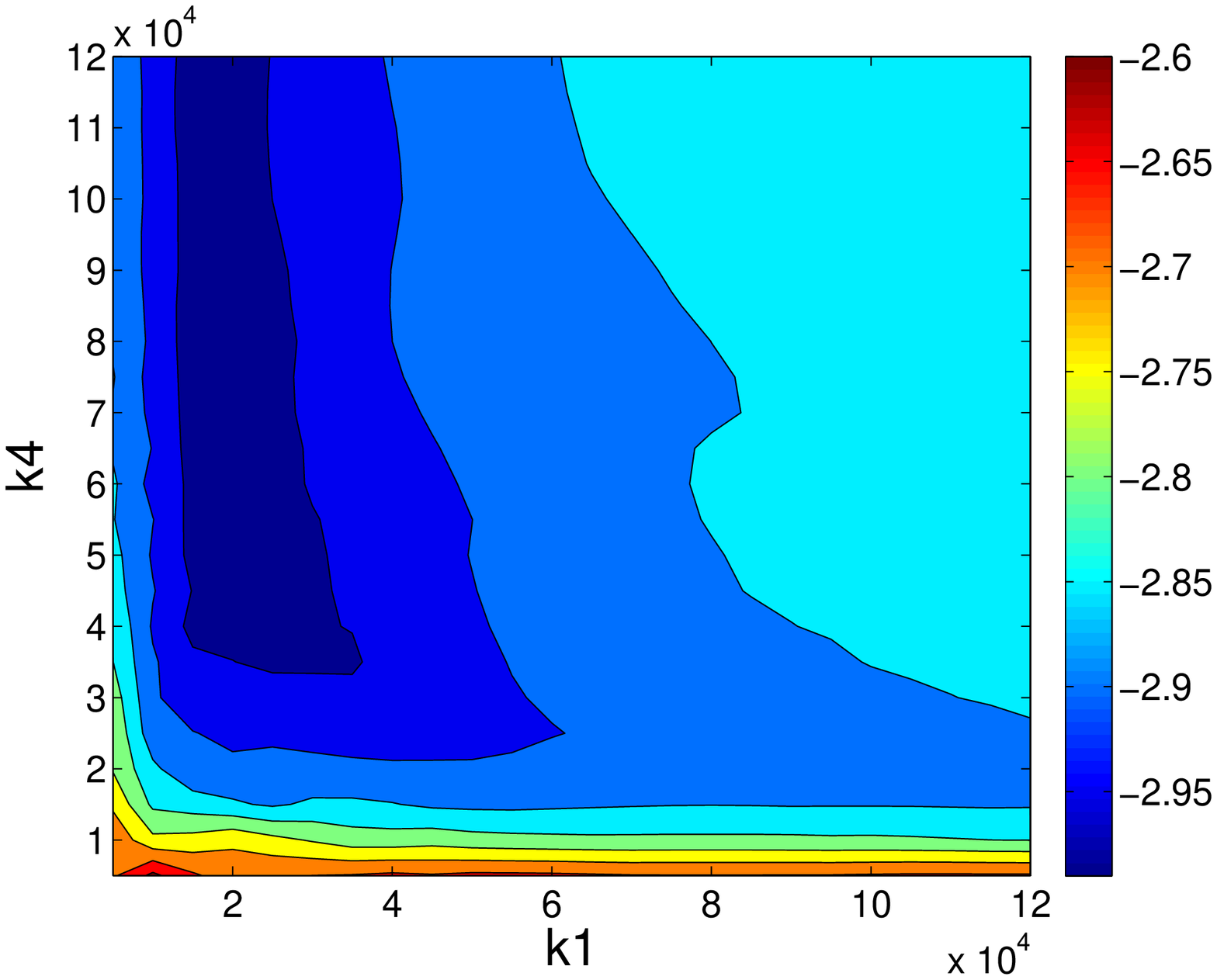}
\includegraphics[width=5.5cm]{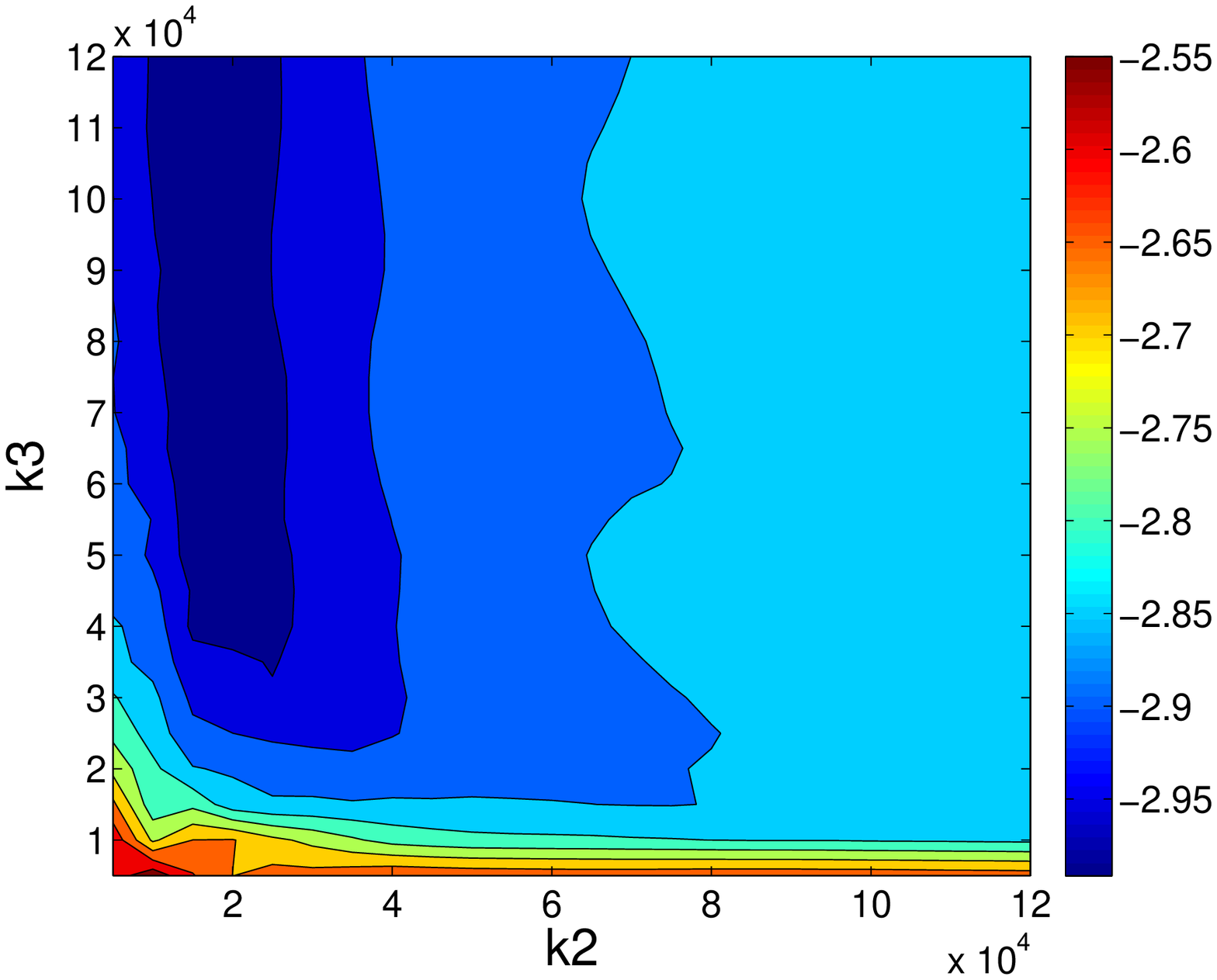}
\includegraphics[width=5.5cm]{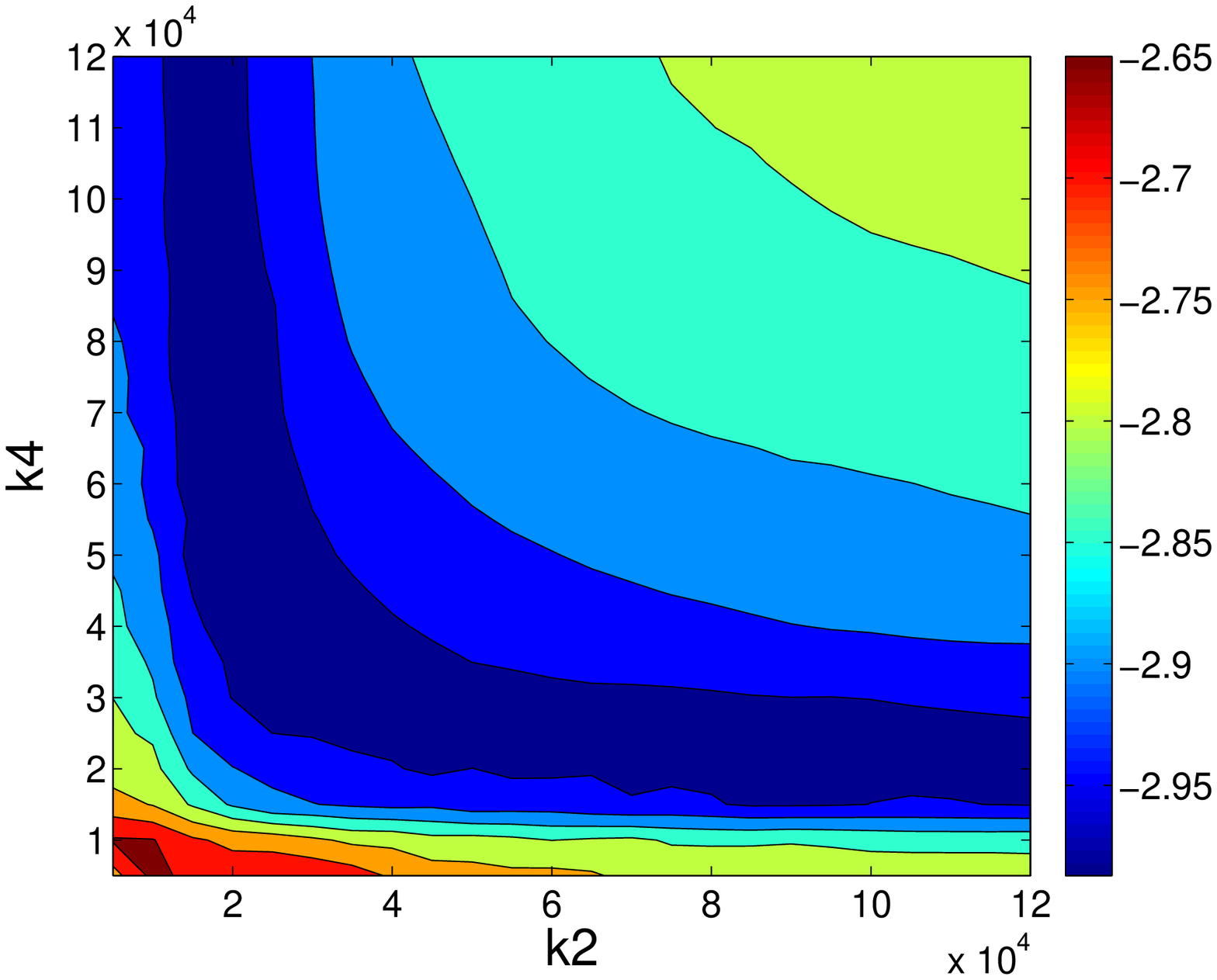}
\includegraphics[width=5.5cm]{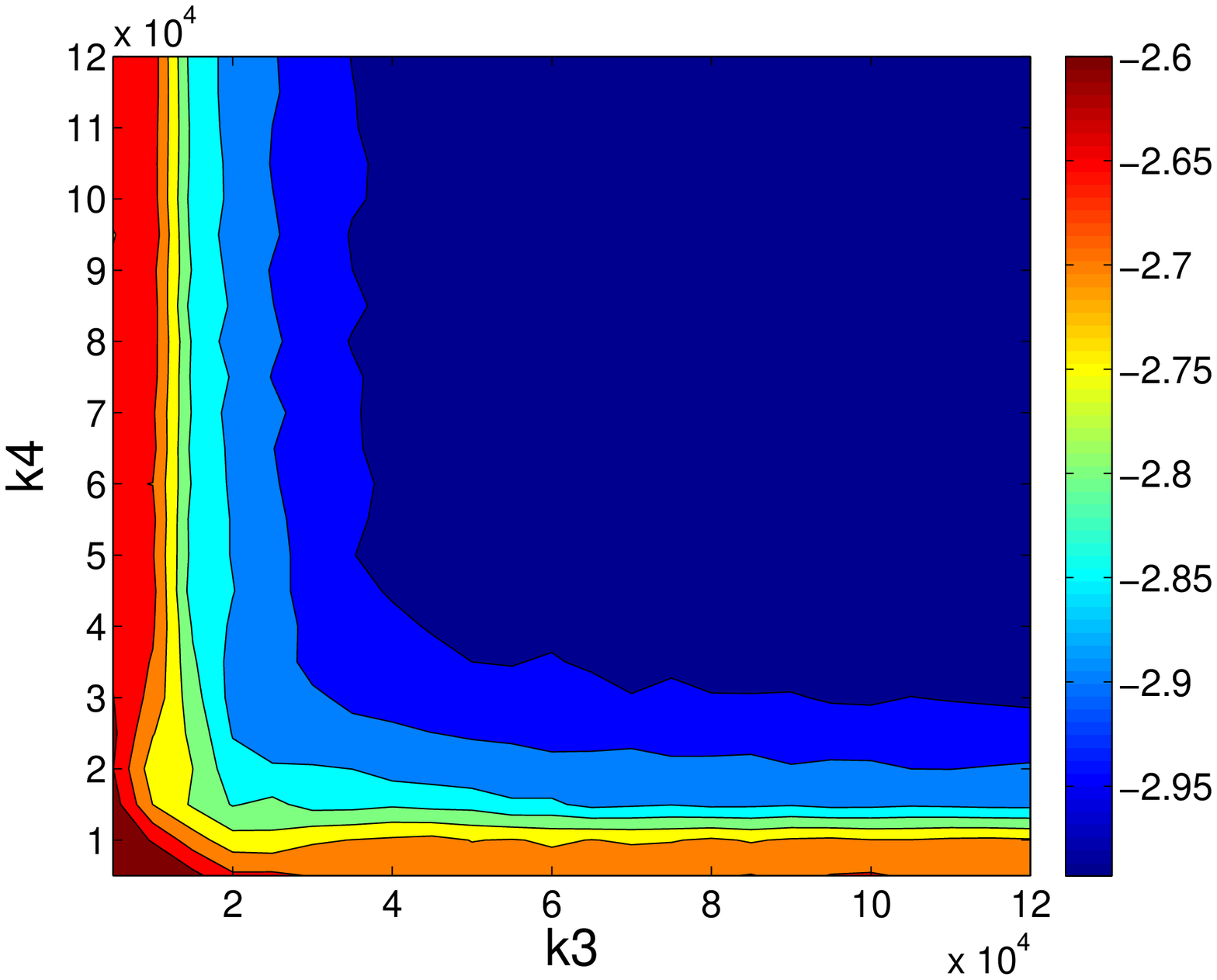}
\caption{Bivariate variations of the exact logarithm of cost function around the optimal value. \label{fig12}}
\end{center}
\end{figure}
\begin{figure}[h!]
\begin{center}
\includegraphics[width=5.5cm]{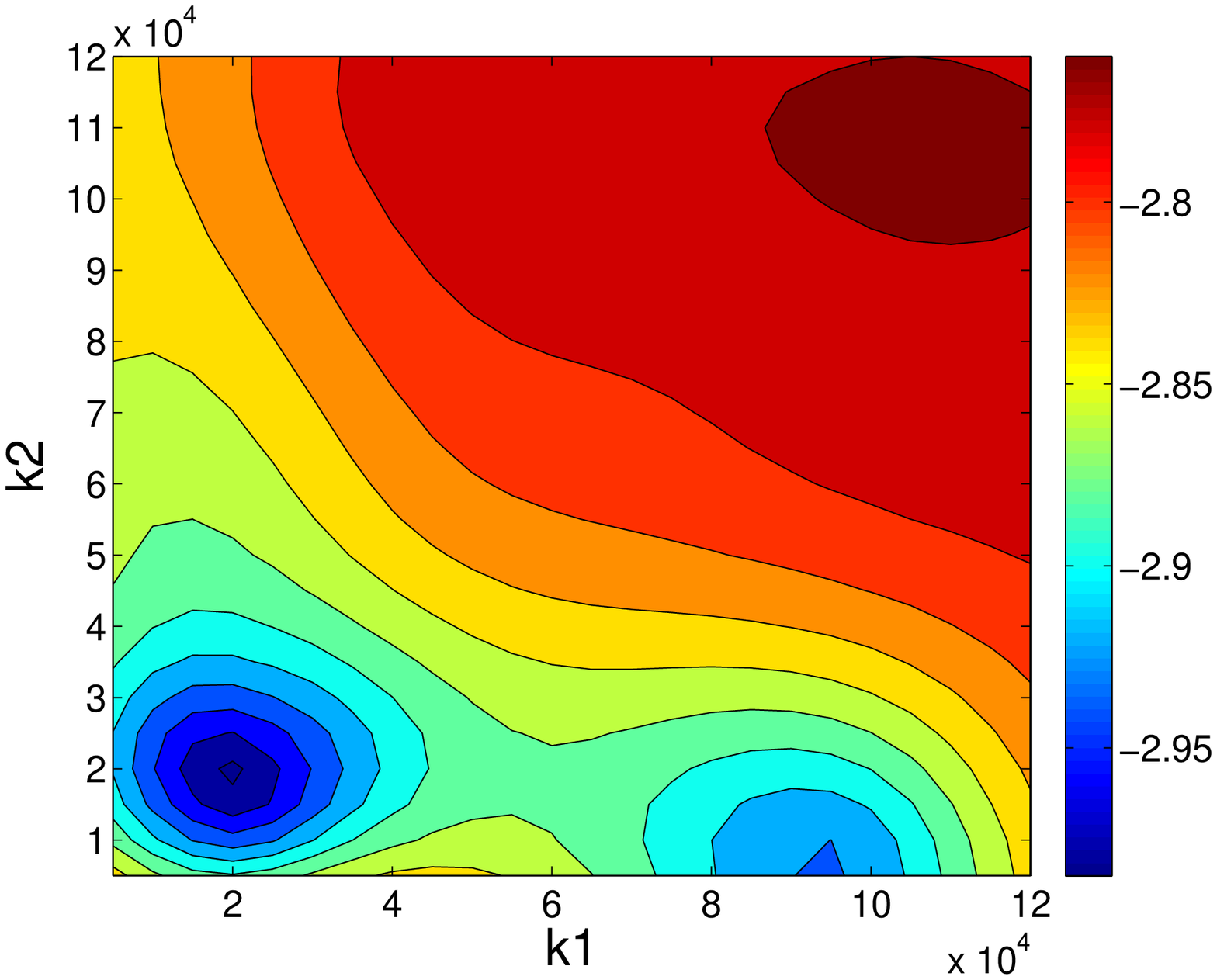}
\includegraphics[width=5.5cm]{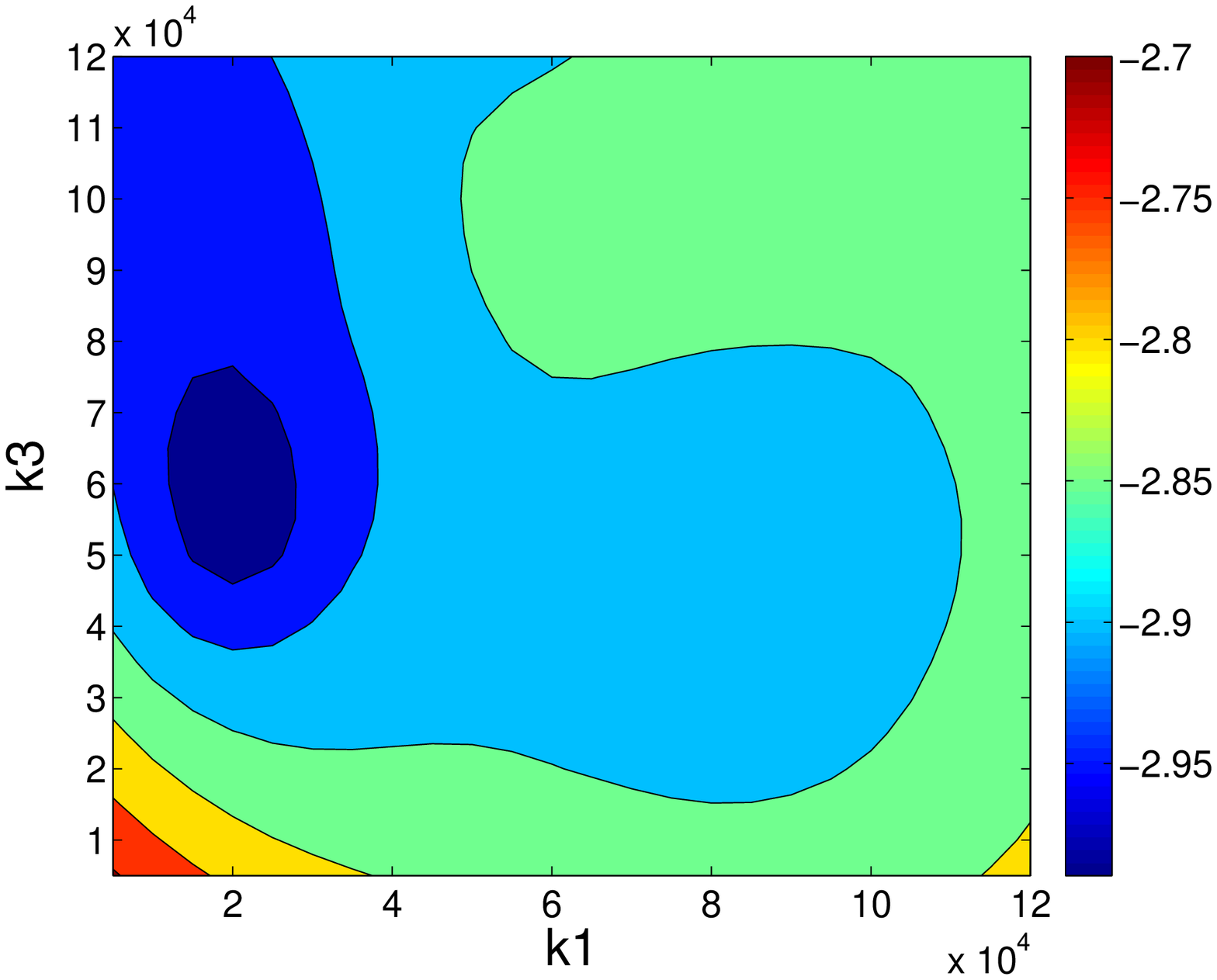}
\includegraphics[width=5.5cm]{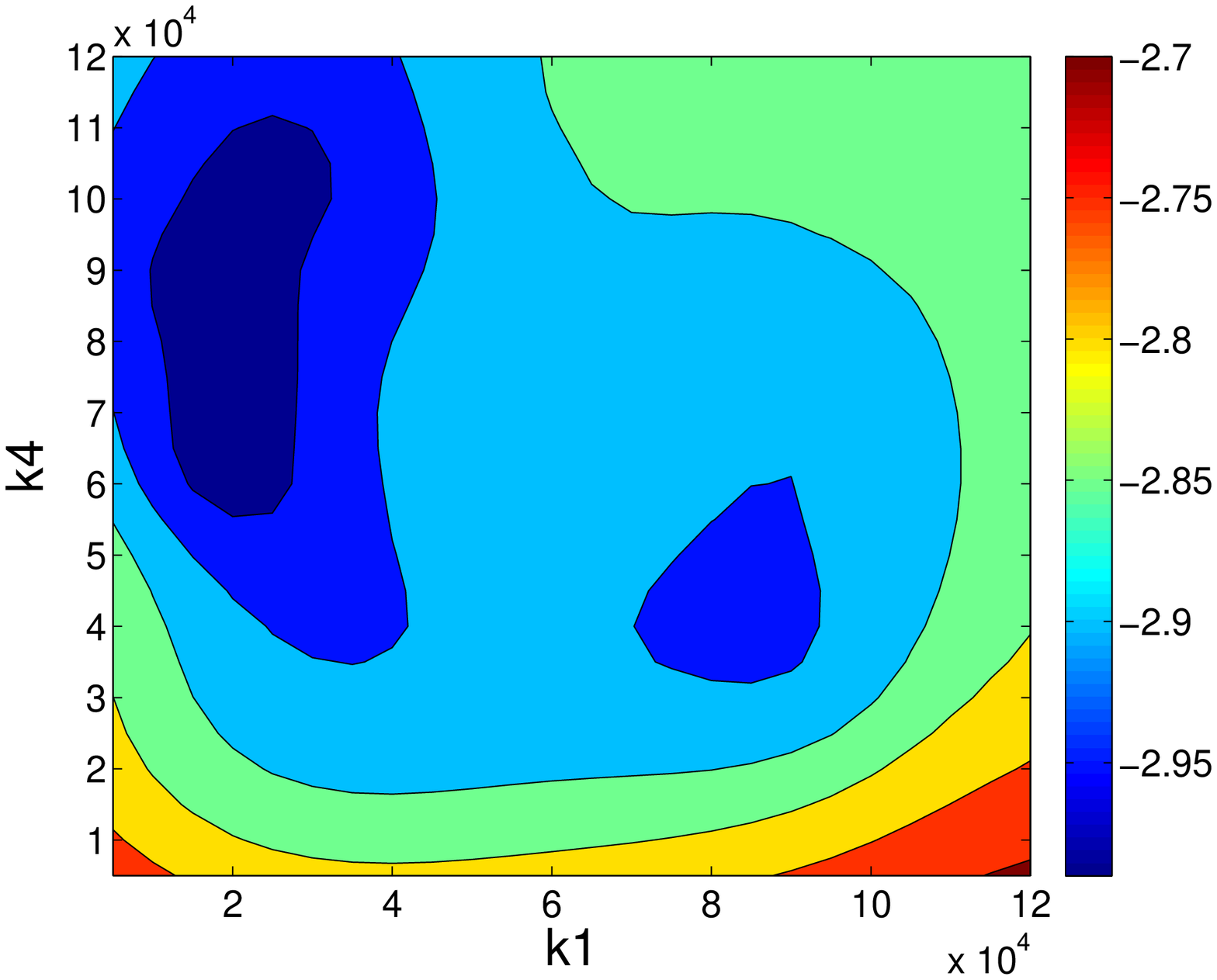}
\includegraphics[width=5.5cm]{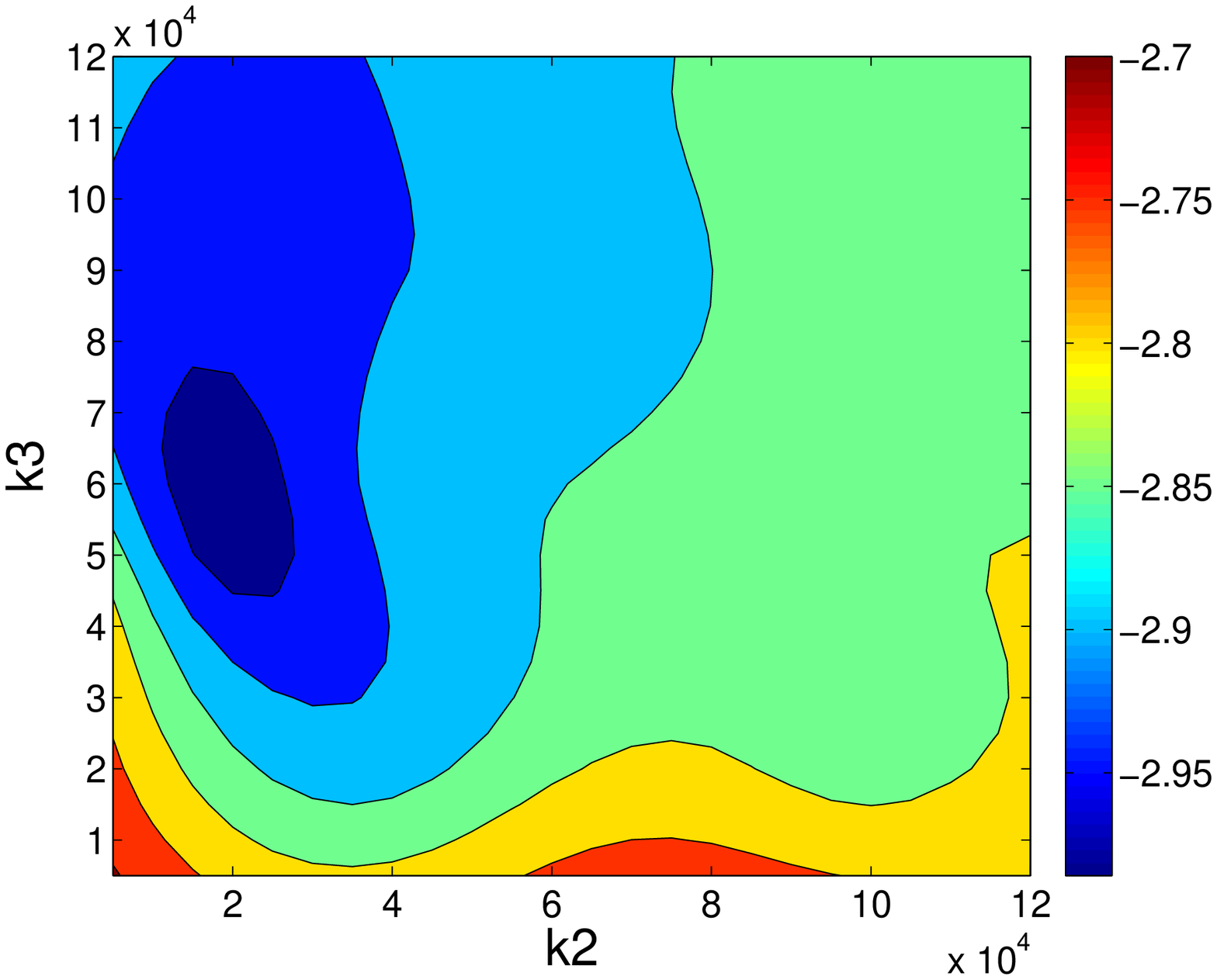}
\includegraphics[width=5.5cm]{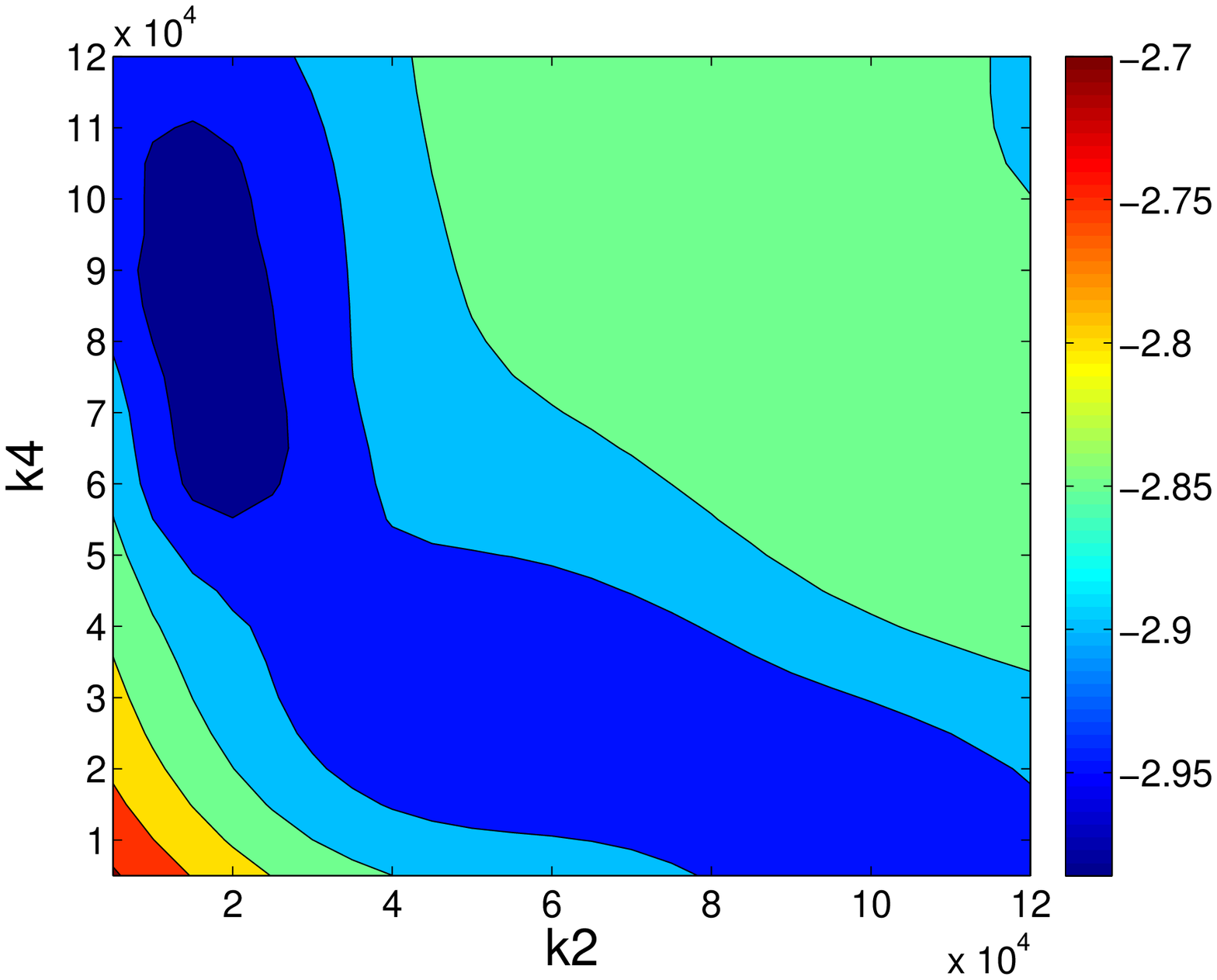}
\includegraphics[width=5.5cm]{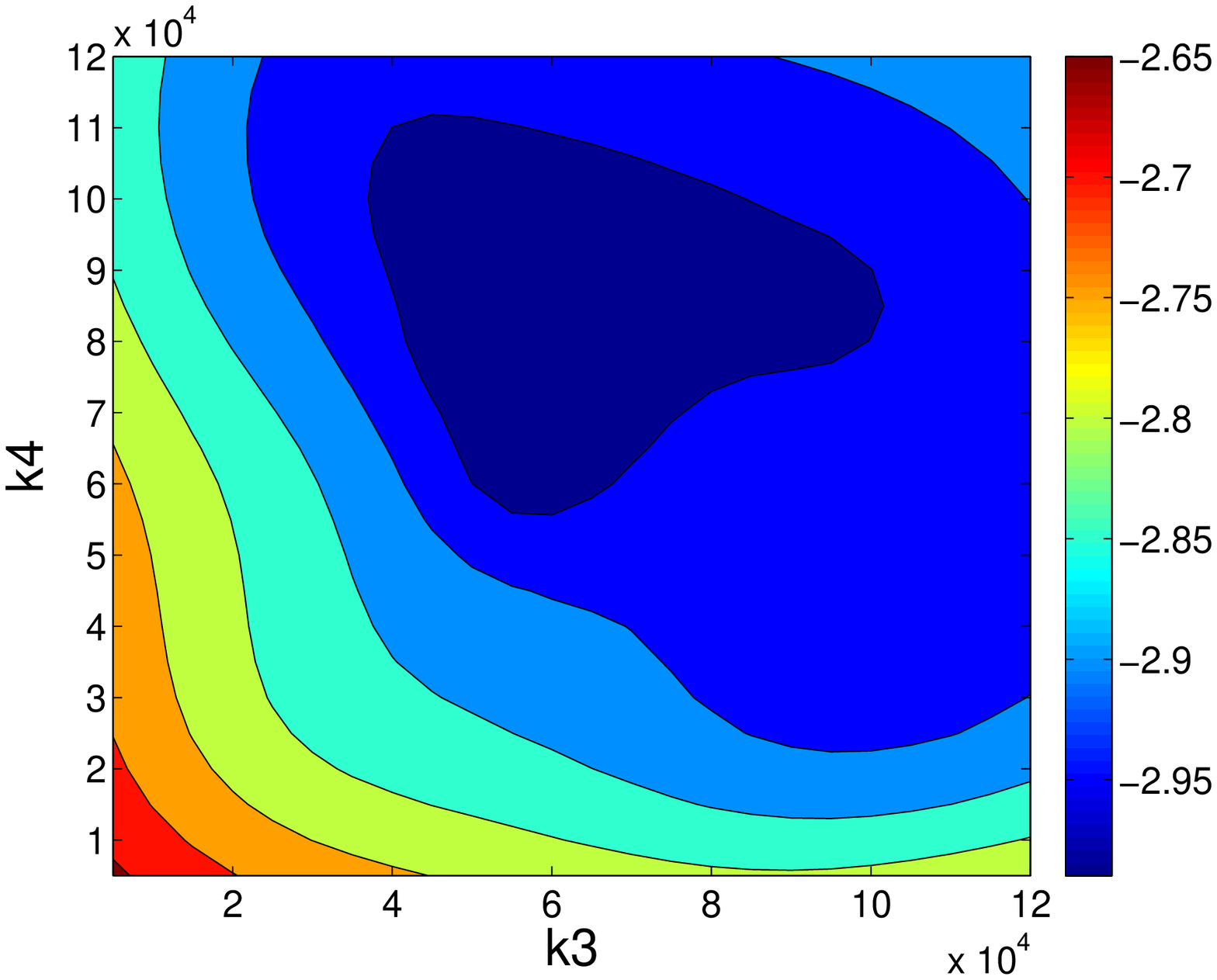}
\caption{Bivariate variations of the approximate logarithm of cost function around the optimal value. \label{fig12b}}
\end{center}
\end{figure}
The new norm of the displacement of the observation point calculated using the designed structure is plotted on Fig.\ref{fig13}. The previous maximum displacement was $5.97\times 10^{-2}$~m. The new one, after optimization, is $4.82\times 10^{-2}$~m.
\begin{figure}[h!]
\begin{center}
\includegraphics[width=8.5cm]{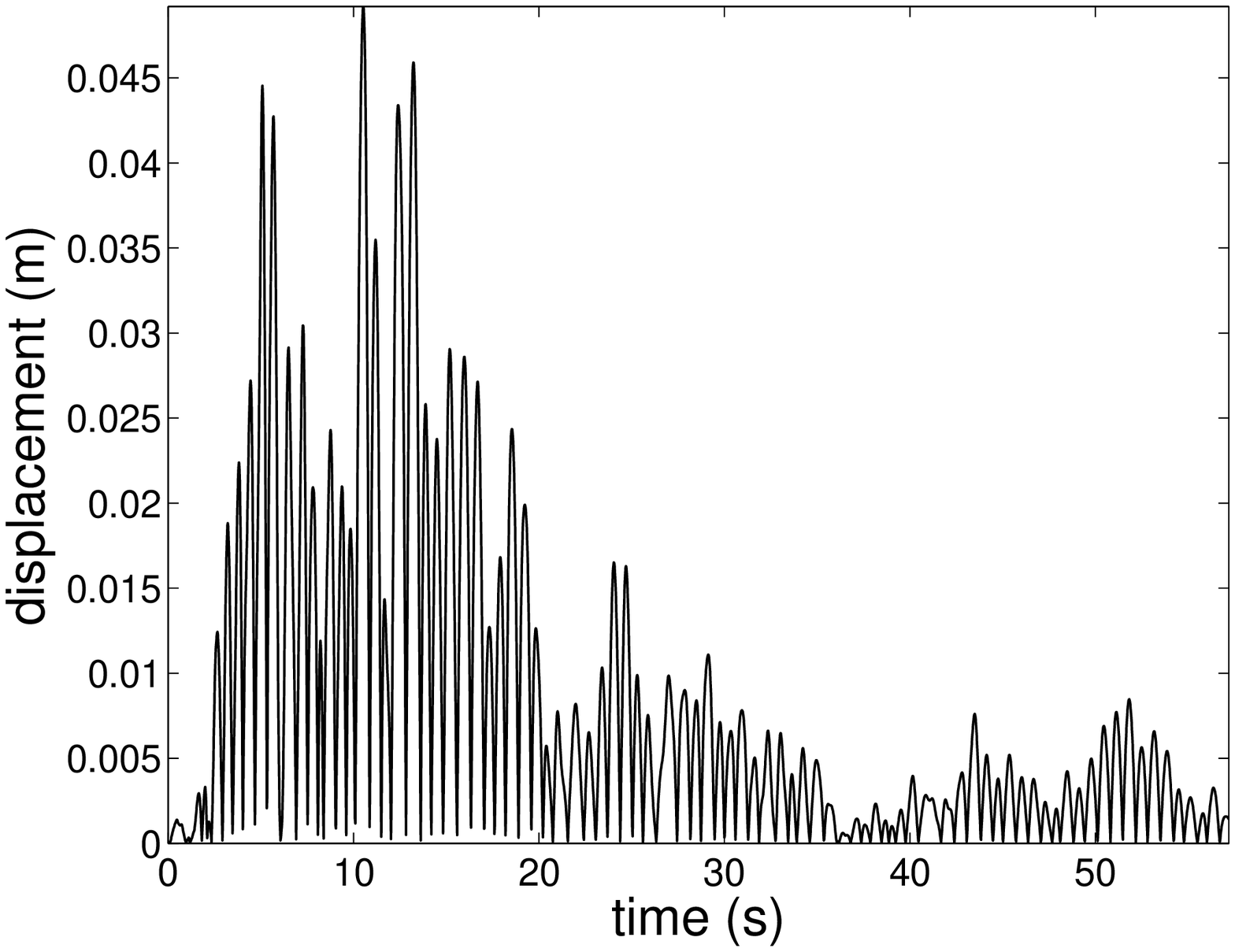}
\caption{Norm of the displacement at observation point for the designed structure. \label{fig13}}
\end{center}
\end{figure}

\section{Conclusions}

We have presented an original variant of the simulated annealing algorithm (1) by incorporating an ISDE generator for the temperature-dependent probability distribution in order to explore efficiently the search space and (2) by interpolating the cost function using polyharmonic splines in order to calculate explicitly and efficiently the cost function and its gradient. Several Markov Chain can be generated in parallel in order to accelerate the enrichment of the interpolation surface. The algorithm has been illustrated on two applications. The first application has shown the potentiality of the algorithm in high dimension and its ability to escape from local minima. The second application shows the potentiality of the algorithm for large computational models.
In addition, all the improvements that have been proposed in \cite{Ingber1993,Pardalos2002} for the classical simulated annealing algorithm could be used to improve the algorithm proposed in the present paper. Two remaining issues will be treated in future works. The first one concerns the choice the temperature law which has not been discussed in this paper since the convergence properties for classical simulated annealing are not preserved exactly. The second issue concerns the choice of the size step for the integration of the ISDE. Since the cost function is know explicitly by its polyharmonic spline representation, an adaptative step size could be studied and implemented in the algorithm.

\end{document}